\documentclass[a4paper,reqno,10pt]{amsart}
\usepackage{amsmath}
\usepackage{graphicx}
\usepackage{amsfonts}
\usepackage{amssymb}
\newtheorem{theorem}{Theorem}[section]
\newtheorem{lemma}[theorem]{Lemma}
\newtheorem{proposition}[theorem]{Proposition}
\newtheorem{corollary}[theorem]{Corollary}
\newtheorem{hypothesis}[theorem]{Hypothesis}
\theoremstyle{definition}
\newtheorem{definition}[theorem]{Definition}

\theoremstyle{remark}
\newtheorem{remark}[theorem]{Remark}
\numberwithin{equation}{section}
\def\Dim{\noindent\hbox{{\bf Proof.}$\;\; $}}
\def\finedim{{\hfill\hbox{\enspace${ \square}$}} \smallskip}
\def\sqr#1#2{{\vcenter{\vbox{\hrule height .#2pt
\hbox{\vrule width .#2pt height#1pt \kern#1pt \vrule
width .#2pt} \hrule height .#2pt}}}}
\def\square{\mathchoice\sqr54\sqr54\sqr{4.1}3\sqr{3.5}3}
\newcommand{\abs}[1]{\lvert#1\rvert}

\def\nat{\mathbb N}

\def\Pr{\mathbb P}
\def\E{\mathbb E}

\def\R{\reali}

\def\reali{\mathbb R}

\def\F{{\mathcal F}}

\def\Ft{{\mathcal F}_t}
\def\P{{\mathcal P}}

\def\Et{\mathbb{E}^{\F_t}}

\begin{document}
\title{Infinite Horizon and Ergodic Optimal Quadratic Control for an Affine Equation with Stochastic Coefficients}

\author{Giuseppina Guatteri and Federica Masiero}

\address{Giuseppina Guatteri, Dipartimento di Matematica,
  Politecnico di Milano, piazza Leonardo da Vinci 32, 20133 Milano,
e-mail: giuseppina.guatteri@polimi.it}
\address{Federica Masiero,
Dipartimento di Matematica e Applicazioni, Universit\`{a} di Milano
Bicocca, via R. Cozzi 53 - Edificio U5, 20125 Milano, e-mail:
federica.masiero@unimib.it}
\maketitle
\begin{abstract}
We study quadratic optimal stochastic control problems with
control dependent noise state equation perturbed by an affine term
and with stochastic coefficients. Both infinite horizon case and
ergodic case are treated. To this purpose we introduce a Backward
Stochastic Riccati Equation and a {\em dual} backward stochastic
equation, both considered in the whole time line. Besides some
{\em stabilizability} conditions we prove existence of a solution
for the two previous equations defined as limit of suitable finite
horizon approximating problems. This allows to perform the
synthesis of the optimal control.
\end{abstract}



\medskip

\email{}

\bibliographystyle{plain}

%
%
%
%
%
%
%
%
%
%

\section{Introduction}

Backward Stochastic Riccati Equations (BSREs) are naturally linked with stochastic
optimal control problems with stochastic coefficients.
The first existence and uniqueness result for such a kind of equations has been given by Bismut in \cite{Bi76}, but then several works, see \cite{Bi78}, \cite{KoTa01}, \cite{KoTa02}, \cite{KoTa03},
\cite{KoZh}, \cite{Peng} and \cite{Pe99}, followed as the problem, in its general formulation, turned out to be difficult to handle and challenging. Indeed only very recently Tang in \cite{Tang} solved the general non singular case corresponding to the linear quadratic problem with random coefficients and control dependent noise.

In his  paper the so-called linear quadratic optimal control problem is considered: minimize over all admissible controls $u$ the following cost functional
\begin{equation}
J(0,x,u)=\mathbb{E}\int_{0}^{T}[\left\langle S_{s}X_{s}%
,X_{s}\right\rangle +|u_{s}|^{2}]ds + \mathbb{E} \langle PX_T,X_T \rangle ;
\label{costo.fin.intro}
\end{equation}
and $\inf_{u}J(0,x,u)$ is the so called value function.
In (\ref{costo.fin.intro}), $X_s \in \R^n$ is solution of the following linear stochastic system:
\begin{equation}
\left\{
\begin{array}
[c]{ll}
dX_{s}=(A_{s}X_{s}+B_{s}u_{s})ds+
{\displaystyle\sum_{i=1}^{d}}
\left(  C_{s}^{i}X_{s}+D_{s}^{i}u_{s}\right)  dW_{s}^{i} & s\geq 0\\
X_{0}=x, &
\end{array}
\right.  \label{stato.Tang.intro}%
\end{equation}
where $W$ is a $d$ dimensional brownian motion and $A, B, C, D, S$
are stochastic processes predictable with respect to the natural augmented filtration
$\{ \F_t\}_{t \geq 0}$ while $P$ is a random variable
$\F_T$-measurable.

All these results cover the finite horizon case.
\medskip

In this paper starting from the results of \cite{Tang}, we address the infinite horizon case and the ergodic case.
The infinite horizon case, with random coefficients, has been addressed also in \cite{GTinf}.
In \cite{GTinf} the infinite dimensional case is treated but no control dependent noise appears in the state equation.  Since our final goal is to address ergodic control, in the state equation we consider a forcing term. Namely, the state equation that describe the system under control is the following affine stochastic equation:
\begin{equation}
\left\{
\begin{array}
[c]{ll}
dX_{s}=(A_{s}X_{s}+B_{s}u_{s})ds+
{\displaystyle\sum_{i=1}^{d}}
\left(  C_{s}^{i}X_{s}+D_{s}^{i}u_{s}\right)  dW_{s}^{i} +f_s ds & s\geq 0\\
X_{0}=x, &
\end{array}
\right.  \label{f.stato.stato.intro}%
\end{equation}
Our main goal  is to minimize with respect to $u$ the infinite
horizon cost functional,
\begin{equation}
J_{\infty}(0,x,u)=\mathbb{E}\int_{0}^{+\infty}[\left\langle S_{s}X_{s}%
,X_{s}\right\rangle +|u_{s}|^{2}]ds
\label{costo.inf.intro}
\end{equation}
and the following ergodic cost functional:
\begin{equation}\liminf_{\alpha \to 0} \alpha J^{\alpha}(0,x,u)
\label{costo.erg.intro}
\end{equation}
where
\begin{equation}
J^{\alpha}(0,x,u)=\mathbb{E}\int_{0}^{+\infty}e^{-2\alpha s}[\left\langle S_{s}X_{s}%
,X_{s}\right\rangle +|u_{s}|^{2}]ds,
\label{costo_scontato.intro}
\end{equation}
In order to carry on this programme we have first to reconsider the
finite horizon case since now the state equation is affine. As it
is well known the value function has in the present situation a
quadratic term represented in term of the solution of the Backward
Stochastic Riccati Equation (BSRE) in $[0,T]$:
\begin{align}
&dP_{t}    =-\left[  A_{t}^{\ast}P_{t}+P_{t}A_{t}+S_{t}+
{\displaystyle\sum_{i=1}^{d}} \left(  \left(  C_{t}^{i}\right)
^{\ast}P_{t}C_{t}^{i}+\left(  C_{t} ^{i}\right)
^{\ast}Q_{t}+Q_{t}C_{t}^{i}\right)  \right]  dt+
{\displaystyle\sum_{i=1}^{d}}
Q_{t}^{i}dW_{t}^{i}+\label{Riccati0T.intro}, \\
&  \left[  P_tB_t+ {\displaystyle\sum_{i=1}^{d}} \left(  \left(
C^{i}_t\right)  ^{\ast}P_t D^{i}_t +Q^{i}D_t ^{i}\right)  \right]
\left[  I+ {\displaystyle\sum_{i=1}^{d}} \left(  D^{i}_t\right)
^{\ast}P_tD^{i}_t\right]^{-1}\!\!\left[  P_tB_t+
{\displaystyle\sum_{i=1}^{d}} \left(  \left(  C^{i}_t\right)
^{\ast}P_t D_t ^{i}+Q^{i}_t D^{i}_t\right)
\right] ^{\ast}\!\!\,\!\! dt, \nonumber \\
&P_T=P,\nonumber
\end{align}
and a linear term involving the so-called {\em costate equation}
(dual equation):
 \begin{equation}
\left\lbrace
\begin{array}[c]{ll}
dr^T_{t}=-H_{t}^{*}r^T_tdt-P_{t}f_{t}dt-{\displaystyle
\sum_{i=1}^{d}\left( K_{t}^{i} \right)^{*}g_{t}^{i,T}}dt+{\displaystyle
\sum_{i=1}^{d}}g_{t}^{i,T}dW_{t}^{i}, & t \in \left[ 0,T\right]  \\
r_{T}^T=0. &
\end{array}
\right. \label{dualeT.intro}
\end{equation}
The coefficients $H$ and $K$ are related with  the coefficients of
the state equation and the solution to the BSRE on $[0,T]$. In
details, if we denote for $t \in \left[ 0,T \right]$
\begin{align*}
& \Lambda\left( t, P_{t}, Q_{t} \right)=-[
I+{\displaystyle\sum_{i=1}^{d}} \left(  D_{t}^{i}\right)
^{\ast}P_{t}D_{t}^{i}]
^{-1}[P_{t}B_{t}+{\displaystyle\sum_{i=1}^{d}} (
Q_{t}^{i}D_{t}^{i}+\left(  C_{t}^{i}\right)
^{\ast}P_{t}D_{t}^{i})],  ^{\ast}
\end{align*}
then we have: $H_{t}=A_{t}+B_{t}\Lambda\left( t, P_{t}, Q_{t} \right)$
and  $K_{t}^{i}=C_{t}^{i}+D_{t}^{i}\Lambda\left( t, P_{t}, Q_{t}
\right)$.
 The solution $(r,g)$ of this equation together with the
solution $(P,Q)$ of the BSRE equation \eqref{Riccati0T.intro}
allow to describe the optimal control and perform the synthesis of
the optimal equation. Equation \eqref{dualeT.intro} is the
generalization of the deterministic equation considered by
Bensoussan in \cite{Ben} and by Da Prato and Ichikawa in
\cite{DaPIch} and of the stochastic backward equation introduced
in \cite{Tess2} for the case without control dependent noise and
with deterministic coefficients.

The main difference from the equation considered in \cite{Tess2}
is that, being the solution to the Riccati equation a couple of
stochastic processes $(P,Q)$ with $Q$  just square integrable,
equation \eqref{dualeT.intro} has stochastic coefficients that are
not uniformly bounded. So the usual technique of resolution does
not apply directly. When $r$ is one dimensional similar difficulties are 
treated e.g. in \cite{BriCon} using Girsanov
Theorem and properties of BMO martingales. Here being the problem
naturally multidimensional we can not apply the Girsanov
transformation to get rid of the term $\sum_{i=1}^{d}\left(
K_{t}^{i} \right)^{*}g_{t}^{i}dt+{
\sum_{i=1}^{d}}g_{t}^{i}dW_{t}^{i}$.

Nevertheless we can exploit a duality relation between the dual
equation  \eqref{dualeT.intro} and the following equation
\begin{equation}
\left\{
\begin{array}
[c]{ll}
dX_{s}=H_{s}X_sds+{\displaystyle\sum_{i=1}^{d}}K_{s}^{i}X_{s}dW_{s}^{i} & s\in \left[ t,T\right) \\
X_{t}=x.&
\end{array}
\right.  \label{closedloopHK.intro}%
\end{equation}
This equation is indeed the closed loop equation
 related to the linear quadratic problem
and can be solved following \cite{Gal} and its control
interpretation allows to gain enough regularity to perform the
duality relation with $(r,g)$.

Once we are able to handle the finite horizon case, we can proceed
to study the infinite horizon  problem. The BSRE corresponding to
this problem is, for $t \geq 0$,
\begin{align}
&dP_{t}    =-\left[  A_{t}^{\ast}P_{t}+P_{t}A_{t}+S_{t}+
{\displaystyle\sum_{i=1}^{d}}
\left(  \left(  C_{t}^{i}\right)  ^{\ast}P_{t}C_{t}^{i}+\left(  C_{t}
^{i}\right)  ^{\ast}Q_{t}+Q_{t}C_{t}^{i}\right)  \right]  dt+
{\displaystyle\sum_{i=1}^{d}}
Q_{t}^{i}dW_{t}^{i}+\label{RiccatiAlg.intro}\\
&  \left[  P_tB_t+ {\displaystyle\sum_{i=1}^{d}} \left(  \left(
C^{i}_t\right)  ^{\ast}P_t D^{i}_t +Q^{i}D_t ^{i}\right)
\right]\!\! \left[  I+ {\displaystyle\sum_{i=1}^{d}} \left(
D^{i}_t\right) ^{\ast}P_tD^{i}_t\right]  ^{-1}\!\!\left[  P_tB_t+
{\displaystyle\sum_{i=1}^{d}} \left(  \left(  C^{i}_t\right)
^{\ast}P_t D_t ^{i}+Q^{i}_t D^{i}_t\right) \right]
^{\ast}\,\!\!\!\!\! dt.\nonumber
\end{align}
Note that differently from equation \eqref{Riccati0T.intro}, the
final condition has disappeared since the horizon is infinite. We assume 
a suitable stabilizability condition, see also
\cite{GTinf}, namely we ask that there 
exists a control
$u\in L_{\mathcal{P}}^{2}([0,+\infty)\times\Omega;\R^k)$ such that for all
$t\geq0$ and all $x\in\mathbb{R}^{n}$%
\begin{equation*}
\mathbb{E}^{\mathcal{F}_{t}}\int_{t}^{+\infty}[\left\langle S_{s}X_{s}%
^{t,x,u},X_{s}^{t,x,u}\right\rangle +|u_{s}|^{2}]ds<M_{t,x}.
\end{equation*}
for some positive constant $M_{t,x}$. By $L_{\mathcal{P}}^{2}([0,+\infty)\times\Omega;\R^k)$ we denote the space of predictable square integrable processes. 
Under this stabilizability condition, we prove 
that there exists a minimal solution
$(\overline{P},\overline{Q})$ and we can perform the synthesis of
the optimal control with $f=0$. More precisely we introduce a
sequence $(P^N,Q^N)$ of solutions of the Riccati equation in
$[0,N]$ with $P^N(N)=0$ and we show that for any $t \geq 0$ the
sequence of $P^N$ pointwise converge, as $N$ tends to $+\infty$,
to a limit denoted by $\overline{P}$. The sequence of $Q^N$ instead
only converge weakly in $L^2_\P(\Omega\times [0,T])$
 to some process $\overline{Q}$ and this is not enough to pass to the limit in
the fundamental relation and then to conclude that the limit
$(\overline{P}, \overline{Q})$ is the solution for the infinite horizon
Riccati equation \eqref{RiccatiAlg.intro}. Therefore, as for the
finite horizon case, we have to introduce the stochastic
Hamiltonian system to prove that the limit $(\overline{P},\overline{Q})$ solves
the BRSE \eqref{RiccatiAlg.intro}, see Corollary \ref{CorMin}.
Indeed studying the stochastic Hamiltonian system we can prove
that the optimal cost for the approximating problem converge to
the optimal cost of the limit problem and this implies that
$\overline{P}$ is the solution of the BSRE.

In order to cope with the affine term $f$ we have to introduce an
infinite horizon, this time, backward equation
\begin{equation}
dr_{t}=-H_{t}^{*}r_t dt-P_{t}f_{t}dt-{\displaystyle\sum_{i=1}^{d}\left( K_{t}^{i} \right)^{*}g_{t}^{i}}dt+{\displaystyle\sum_{i=1}^{d}}g_{t}^{i}dW_{t}^{i},\text{
\ \ \ \ \ } t \geq 0.  \\
\label{duale.intro}
\end{equation}
Notice that the typical monotonicity assumptions on the
coefficients of this infinite horizon BSDE are replaced by the
finite cost condition and the Theorem of Datko.
As a consequence of this new hypothesis we have that the solution
to the closed loop equation considered in the whole positive time
line with the coefficients evaluated in $\overline{P}$ and $\overline{Q}$,
is exponentially stable.

Hence a solution $(\bar{r},\bar{g})$ to this equation is
obtained as limit of the sequence $(r^T,g^T)$ defined in
\eqref{dualeT.intro}, indeed using duality and the exponential
stability property of the solution to \eqref{closedloopHK.intro},
we can prove that the sequence of $r^T$ and its limit
$\bar{r}$ are uniformly bounded. Hence, having both
$(\overline{P},\overline{Q})$ and $(\bar{r},\bar{g})$,
we can express the optimal control and the value function.

Eventually we come up with the ergodic case: first of all we set
$X^\alpha_s:= e^{-\alpha s}X_s$ and $u^\alpha_s:= e^{-\alpha s}
u_s$ and we notice that the functional $J^{\alpha}(0,x,u)$ can be
written as an infinite horizon functional in terms of $X^\alpha$
and $u^\alpha$:
\[J^{\alpha}(0,x,u)= \mathbb{E}\int_{0}^{+\infty}[\left\langle S_{s}X^\alpha_{s}%
,X^\alpha_{s}\right\rangle +|u^\alpha_{s}|^{2}]ds. \] This allows
us to adapt the previous results on the infinite horizon when
$\alpha >0$ is fixed. \\ Then, in order to study the limit
\eqref{costo.erg.intro}, we need to investigate the behaviour of
$X^\alpha$, of the solution $P^\alpha$ of the Riccati equation
corresponding to
$\overline{J}^{\alpha}(x):=\inf_{u}J^\alpha(0,x,u)$ and the
solutions $(r^\alpha,g^\alpha)$ of the dual equations
corresponding to $H^\alpha, K^\alpha$ and $f^\alpha_t= e^{-\alpha
t}f_t$. In the general case it turns out that the ergodic limit
has the following form:
\begin{multline*}
 \underline{\lim}_{\alpha\rightarrow 0}\alpha \overline{J}^{\alpha}(x)=  \underline{\lim}_{\alpha\rightarrow 0}2\alpha \E\int_0^{+\infty}\langle r^{\alpha}_s, f^{\alpha}_s\rangle ds \\ -\overline{\lim}_{\alpha\rightarrow 0}\alpha\E\int_0 ^{+\infty} \vert (  I+
{\displaystyle\sum_{i=1}^{d}}
\left(  D_{s}^{i}\right)  ^{\ast}P^{\alpha}_{s}D_{s}^{i})  ^{-1}( B_{s}^{*}r^{\alpha}_{s}+{\displaystyle\sum_{i=1}^{d}}\left( D_{s}^{i}\right) ^{*}g_{s}^{\alpha,i})\vert ^2 ds.
\end{multline*}
When the coefficients of the state equation are deterministic
similar problems have already been treated: we cite \cite{Ben1},
\cite{Tess2} and bibliography therein. In \cite{Ben1} in the state
equation all the coefficients are deterministic and no control
dependent noise is studied, while in \cite{Tess2} only the forcing
term $f$ is allowed to be random.


Finally we describe the content of each section: in section 2,
after recalling some results of \cite{Tang}, we solve the finite
horizon case when the state equation is affine: the key point is the solution of the dual equation \eqref{dualeT.intro}, which is studied in paragraph 2.2; in section 3
we solve the infinite horizon case with $f=0$, in section 4 we
study the infinite horizon equation \eqref{duale.intro}, in section 5 we
complete the general infinite horizon case, finally in section 6
we study the ergodic case.

 {\bf Acknowledgments.} The authors wish to thank Philippe Briand
 for the very useful discussions on the role of BMO martingales in the theory of Backward Stochastic
 Differential Equations and Gianmario Tessitore for very useful discussions.
\section{Linear Quadratic optimal control in the finite horizon case
}
Let $(\Omega,\mathcal{E},(\mathcal{F}_t)_{t\geq
0},\mathbb{P})$ be a stochastic base verifying the usual
conditions, and let $(\mathcal{F}_t)_{t\geq
0}$ be natural augmented filtration generated by the Brownian motion. In $(\Omega,\mathcal{E},\mathbb{P})$ we consider the following stochastic differential equation:
\begin{equation}
\left\{
\begin{array}
[c]{ll}
dX_{s}=(A_{s}X_{s}+B_{s}u_{s})ds+
{\displaystyle\sum_{i=1}^{d}}\left(  C_{s}^{i}X_{s}+D_{s}^{i}u_{s}\right)  dW_{s}^{i}+f_{s}ds & s\in \left[ t,T \right] \\
X_{t}=x &
\end{array}
\right.  \label{stato.T}
\end{equation}
where $X$ is a process with values in $\mathbb{R}^{n}$ and
represents the \textit{state} of the system and is our unknown,
$u$ is a process with values in $\mathbb{R}^{k}$ and represents
the \textit{control}, $\left\{ W_{t}:=\left(
W_{t}^{1},...,W_{t}^{d}\right)  ,\text{ }t\geq0\right\}  $ is a
$d$-dimensional standard Brownian motion and the
initial data $x$ belongs to $\mathbb{R}^{n}$. To stress dependence
of the state $X$ on $u$, $t$ and $x$ we will denote the solution
of equation (\ref{stato.T}) by $X^{t,x,u}$ when needed. The norm
and the scalar product in any finite dimensional Euclidean space
$\mathbb{R}^{m}$, $m\geq1$, will be denoted respectively by
$\left| \cdot\right|  $ and $\left\langle \cdot,\cdot\right\rangle
$.

Our purpose is to minimize with respect to $u$ the cost functional,
\begin{equation}
J(0,x,u)=\E \int_{0}^{T}[\left\langle S_{s}X_{s}%
^{0,x,u},X_{s}^{0,x,u}\right\rangle +|u_{s}|^{2}]ds +\E \mid X_{T}^{0,x,u}\mid ^{2}.  \label{costo.T}
\end{equation}
We also introduce the following random variables, for $t\in\lbrack0,T]$:%
\begin{equation}
J(t,x,u):=\Et \int_{t}^{T}[\left\langle S_{s}X_{s}%
^{0,x,u},X_{s}^{0,x,u}\right\rangle +|u_{s}|^{2}]ds +\Et \mid X_{T}^{0,x,u}\mid ^{2}.  \label{costo.T.t}
\end{equation}

We make the following assumptions on $A$, $B$, $C$ and $D$.
\begin{hypothesis}
\label{genhypT} $\ $

\begin{enumerate}
\item [ A1)]$A:\left[  0,T\right]   \times\Omega\rightarrow$
$\mathbb{R}^{n\times n}$, $B:\left[  0,T\right]
\times\Omega \rightarrow$ $\mathbb{R}^{n\times k}$, $C^{i}:\left[
0,T\right] \times\Omega\rightarrow$ $\mathbb{R}^{n\times
n}$, $i=1,...,d$ and $D^{i}:\left[  0,T\right]
\times\Omega\rightarrow$ $\mathbb{R}^{n\times k}$, $i=1,...,d$,
are uniformly bounded processes predictable.

\item[ A2)] $S:\left[  0,T\right]  \times\Omega\rightarrow
\mathbb{R}^{n\times n}$ is uniformly bounded and predictable and it is
almost surely and almost everywhere symmetric and nonnegative.

\item[ A3)] $f:\left[  0,T\right]  \times\Omega \rightarrow
\mathbb{R}^{n}$ is predictable and $f \in 
L^{\infty} \left(   \left[  0,T\right]  \times \Omega  \right)$. We denote such space of processes $L^{\infty}_\P \left(   \left[  0,T\right]  \times \Omega  \right)$.

\end{enumerate}
\end{hypothesis}
\subsection{Preliminary results on the unforced case}

Next we recall some results obtained in \cite{Tang} for the finite
horizon case, with $f=0$ in equation (\ref{stato.T}). In that paper a finite horizon control problem was
studied, namely minimize the quadratic cost functional
\[
J(0,x,u)=\mathbb{E}\langle PX_{T}^{0,x,u},X_{T}^{0,x,u}\rangle
+\mathbb{E}\int_{0}^{T}[\left\langle S_{s}X_{s}^{0,x,u},X_{s}^{0,x,u}%
\right\rangle +|u_{s}|^{2}]ds,
\]
where $P$ is a random matrix uniformly bounded and almost surely positive and
symmetric, $T>0$ is fixed and $X^{0,x,u}$ is the solution to equation (\ref{stato.T}) with $f=0$. To this controlled problem, the following
(finite horizon) backward stochastic Riccati differential equation (BSRDE\ in
the following) is related:
\begin{equation}
\left\{
\begin{array}
[c]{l}
-dP_{t}=G\left(  A_{t},B_{t},C_{t},D_{t};S_{t};P_{t},Q_{t}\right)  dt+
{\displaystyle\sum_{i=1}^{d}}

Q_{t}^{i}dW_{t}^{i}\\
P_{T}=P.
\end{array}
\right.  \label{Riccati0T}%
\end{equation}
where
\[
G\left(  A,B,C,D;S;P,Q\right)  =A^{\ast}P+PA+S+\displaystyle\sum_{i=1}^{d}\left(  \left(  C^{i}\right)  ^{\ast}PC^{i}+\left(  C^{i}\right)  ^{\ast}Q+QC^{i}\right)  -G_{1}\left(  B,C,D;P,Q\right)  ,
\]
and
\begin{align*}
G_{1}\left(  B,C,D;P,Q\right)   &  =\left[  PB+{\displaystyle\sum_{i=1}^{d}}
\left(  \left(  C^{i}\right)  ^{\ast}PD^{i}+Q^{i}D^{i}\right)  \right]
\times\left[  I+{\displaystyle\sum_{i=1}^{d}}
\left(  D^{i}\right)  ^{\ast}PD^{i}\right]  ^{-1}\times\\
&  \times\left[  PB+
{\displaystyle\sum_{i=1}^{d}}
\left(  \left(  C^{i}\right)  ^{\ast}PD^{i}+Q^{i}D^{i}\right)  \right]
^{\ast}
\end{align*}

\begin{definition}
\label{defgen}A pair of predictable processes $\left(  P,Q\right)  $ is a\emph{
solution} of equation (\ref{Riccati0T}) if

\begin{enumerate}
\item $
{\displaystyle\int_{0}^{T}}
\left|  Q_{s}\right|  ^{2}ds<+\infty$, almost surely,

\item
\[
{\displaystyle\int_{0}^{T}}\left|  G\left(  A_{s},B_{s},C_{s},D_{s};S_{s};P_{s},Q_{s}\right)  \right|
ds<+\infty,
\]

\item  for all $t\in\left[  0,T\right]  $
\[
P_{t}=P+{\displaystyle\int_{t}^{T}}G\left(  A_{s},B_{s},C_{s},D_{s};S_{s};P_{s},Q_{s}\right)  ds-{\displaystyle\int_{t}^{T}}{\displaystyle\sum_{i=1}^{d}}Q_{s}^{i}dW_{s}^{i}.%
\]
\end{enumerate}
\end{definition}
\begin{theorem}[\cite{Tang}, Theorems 3.2 and 5.3]
Assume that $A,B,C,D$ and $S$ verify hypothesis \ref{genhypT}.
Then there exists a unique solution to equation (\ref{Riccati0T}).
Moreover the following fundamental relation holds true, for all
$0\leq t\leq s\leq T$, and all $u\in L_{\mathcal{P}}^{2}\left(
\left[ 0,T\right]  \times \Omega,\mathbb{R}^{k}\right)  $:
\begin{align}
\left\langle P_{t}x,x\right\rangle  &  =\mathbb{E}^{\mathcal{F}_{t}%
}\left\langle PX_{T}^{t,x,u},X_{T}^{t,x,u}\right\rangle +\mathbb{E}%
^{\mathcal{F}_{t}}\int_{t}^{T}[\left\langle S_{r}X_{r}^{t,x,u},X_{r}%
^{t,x,u}\right\rangle +|u_{r}|^{2}]dr\label{fundrelT}\\
&  -\mathbb{E}^{\mathcal{F}_{t}}\int_{t}^{T}\left|  \left(  I+{\displaystyle\sum_{i=1}^{d}}\left(  D_{s}^{i}\right)  ^{\ast}P_{s}D_{s}^{i}\right)  ^{1/2}\right.\ast\nonumber\\
&  \ast\left.  \left[  u_{s}+\left(  I+{\displaystyle\sum_{i=1}^{d}}\left(  D_{s}^{i}\right)  ^{\ast}P_{s}D_{s}^{i}\right)  ^{-1}\left(P_{s}B_{s}+{\displaystyle\sum_{i=1}^{d}}
\left(  Q_{s}^{i}D_{s}^{i}+\left(  C_{s}^{i}\right)  ^{\ast}P_{s}D_{s}%
^{i}\right)  \right)  ^{\ast}X_{s}^{t,x,u}\right]  \right|  ^{2}ds\nonumber
\end{align}
Then the value function is given by
\[
\left\langle P_{0}x,x\right\rangle =\inf_{u \in L_{\mathcal{P}}^{2}\left(  \left[
0,T\right]  \times\Omega,\mathbb{R}^{k}\right)  }\mathbb{E}^{\mathcal{F}_{t}%
}\left\langle PX_{T}^{t,x,u},X_{T}^{t,x,u}\right\rangle +\mathbb{E}%
^{\mathcal{F}_{t}}\int_{t}^{T}[\left\langle S_{r}X_{r}^{t,x,u},X_{r}%
^{t,x,u}\right\rangle +|u_{r}|^{2}]dr
\]
and the unique optimal control has the following closed form:%
\[
\overline{u}_{t}=-\left(  I+{\displaystyle\sum_{i=1}^{d}}
\left(  D_{t}^{i}\right)  ^{\ast}P_{t}D_{t}^{i}\right)  ^{-1}\left(
P_{t}B_{t}+{\displaystyle\sum_{i=1}^{d}}\left(  Q_{t}^{i}D_{t}^{i}+\left(  C_{t}^{i}\right)  ^{\ast}P_{t}D_{t}%
^{i}\right)  \right)  ^{\ast}X_{t}^{0,x,\bar{u}}.
\]
If $\overline{X}$ is the solution of the state equation corresponding to
$\overline{u}$ (that is the optimal state), then $\overline{X}$ is the unique
solution to the\emph{ closed loop }equation:
\begin{equation}
\left\{
\begin{array}
[c]{ll}
d\overline{X}_{s}=\left(A_{s}\overline{X}_{s}-B_{s}(  I+{\displaystyle\sum_{i=1}^{d}}\left(  D_{t}^{i}\right)  ^{\ast}P_{t}D_{t}^{i})  ^{-1}(P_{t}B_{t}+{\displaystyle\sum_{i=1}^{d}}
\left(  Q_{t}^{i}D_{t}^{i}+(  C_{t}^{i})  ^{\ast}P_{t}D_{t}
^{i}\right)  )  ^{\ast}\overline{X}_{s}\right)ds +& \\
{\displaystyle\sum_{i=1}^{d}}\left(  C_{s}^{i}\overline{X_{s}}-D_{s}^{i}(  I+{\displaystyle\sum_{i=1}^{d}}\left(  D_{t}^{i}\right)  ^{\ast}P_{t}D_{t}^{i})  ^{-1}(P_{t}B_{t}+{\displaystyle\sum_{i=1}^{d}}
(  Q_{t}^{i}D_{t}^{i}+\left(  C_{t}^{i}\right)  ^{\ast}P_{t}D_{t}^{i})  )  ^{\ast}\overline{X}_{s})\right)  dW_{s}^{i},
\\ \overline{X}_{t}=x
\end{array}
\right.  \label{loop.ban}
\end{equation}
The optimal cost is therefore given in term of the solution of the Riccati matrix \begin{equation}\label{controllo_ottimo_Tang}J(0,x,\overline{u})=\langle P_{0}x,x
\rangle.\end{equation}
and also the following identity holds, for all $ t \in [0,T]$:
\begin{equation}\label{controllo_ottimo_Tang_t}J(t,x,\overline{u})=\langle P_{t}x,x
\rangle.\end{equation}
\end{theorem}
For $t \in \left[ 0,T \right] $, we denote by
\begin{align}
& \Lambda\left( t, P_{t}, Q_{t} \right)=-\left(  I+{\displaystyle\sum_{i=1}^{d}}\left(  D_{t}^{i}\right)  ^{\ast}P_{t}D_{t}^{i}\right)  ^{-1}\left(P_{t}B_{t}+{\displaystyle\sum_{i=1}^{d}}\left(  Q_{t}^{i}D_{t}^{i}+\left(  C_{t}^{i}\right)  ^{\ast}P_{t}D_{t}^{i}\right)  \right)  ^{\ast}, \nonumber\\
& H_{t}=A_{t}+B_{t}\Lambda\left( t, P_{t}, Q_{t} \right), \nonumber\\
&K_{t}^{i}=C_{t}^{i}+D_{t}^{i}\Lambda\left( t, P_{t}, Q_{t} \right). \nonumber \\ 
\label{notazionifHK}
\end{align}
So the closed loop equation \eqref{loop.ban} can be rewritten as
\begin{equation}
\left\{
\begin{array}
[c]{ll}
dX_{s}=H_{s}X_sds+{\displaystyle\sum_{i=1}^{d}}K_{s}^{i}X_{s}dW_{s}^{i} & s\in \left[ t,T\right) \\
X_{t}=x &
\end{array}
\right.  \label{closedloopHK}%
\end{equation}
It is well known, see e.g. \cite{Gal}, that equation \eqref{closedloopHK} admits a solution.
\begin{remark}
\label{stimeHK}
$\Lambda$, $H$ and $K$ defined in (\ref{notazionifHK}) are related to the feedback operator in the solution of the finite horizon optimal control problem with $f=0$. By the boundedness of $P$ and thanks to standard estimates on the Riccati equation, see \cite{Tang}, theorem 5.3, it turns out that for every stopping time $0\leq \tau \leq T$ a.s.,
\begin{equation*}
\E^{\F_{\tau}} \displaystyle\int_{\tau}^{T}\vert \Lambda\left( t, P_{t}, Q_{t} \right) \vert^{2}dt \leq C ,
\end{equation*}
where $C$ is a constant depending on $T$ and $x$. Since $A$, $B$, $C$ and $D$ are bounded, this property holds true also for $H$ and $K$:
\begin{equation*}
\E^{\F_{\tau}} \displaystyle\int_{\tau}^{T}\vert H_{t} \vert^{2}dt+\E^{\F_{\tau}} \displaystyle\int_{\tau}^{T}\vert K_{t} \vert^{2}dt \leq C ,
\end{equation*}
where now $C$ is a constant depending on $T$, $x$, $A$, $B$, $C$ and $D$. In particular, $\Lambda$, $H$ and $K$ are square integrable. In the following we denote by
\begin{equation}
 C_{H}=\sup_{\tau}\E^{\F_{\tau}} \displaystyle\int_{\tau}^{T}\vert H_{t} \vert^{2}dt,\qquad 
C_{K}=\sup_{\tau}\E^{\F_{\tau}} \displaystyle\int_{\tau}^{T}\vert K_{t} \vert^{2}dt. 
\label{BMOnormaHK}
\end{equation}
where the supremum is taken over all stopping times $\tau$, $\tau\in[0,T]$ a.s..
\end{remark}
\subsection{Costate equation and finite horizon affine control}
In order to solve the optimal control problem related to the nonlinear controlled equation \ref{stato.T}, we introduce the so called \textit{dual equation}, or \textit{costate equation},
\begin{equation}
\left\lbrace
\begin{array}[c]{ll}
dr_{t}=-H_{t}^{*}r_tdt-P_{t}f_{t}dt-{\displaystyle\sum_{i=1}^{d}\left( K_{t}^{i} \right)^{*}g_{t}^{i}}dt+{\displaystyle\sum_{i=1}^{d}}g_{t}^{i}dW_{t}^{i}, & t \in \left[ 0,T\right]  \\
r_{T}=0. &
\end{array}
\right. \label{dualeT}
\end{equation}
At a first step we look for a pair of
predictable processes $\left( r,g \right)$ satisfying equation \eqref{dualeT}, and s.t. $r\!\in\!
L^{\infty}_{loc}\left(\Omega,C([0,T],
\R^{n})\right)$ and $g^{i} \in L^{2}_{loc}\left(\left[ 0,T\right]
\times\Omega, \R^{n}\right)$, for $i=i,...,d$.
$L^{\infty}_{loc}\left(\Omega,C([0,T],
\R^{n})\right)$ is the space of predictable processes $r$ with
values in $\R^{n}$ that admit a continuous version and such that
\begin{align*}
\Pr\left(\sup_{t \in[0,T]}\vert r_{t}\vert <\infty\right)=1.
\end{align*}
$L^{2}_{loc}\left(\left[ 0,T\right] \times\Omega, \R^{n}\right)$ is the space of predictable processes $g$ with values in $\R^{n}$ such that
\begin{align*}
\Pr\left( \int_{0}^{T}\vert g_{s}\vert^{2}ds <\infty\right)=1.
\end{align*}

\begin{lemma}
\label{esistenzadualeT} The backward equation \eqref{dualeT}
admits a unique solution $\left( r,g \right)$ that belongs to the
space $L_{loc}^{\infty}\left(\Omega,C([0,T],
\R^{n})\right) \times L^{2}_{loc}\left(\left[ 0,T\right]
\times\Omega, \R^{n \times d} \right)$.

\end{lemma}
\Dim
In order to construct a solution to equation \eqref{dualeT}, we essentially follow \cite{Yo-Z}, chapter 7, where linear BSDEs with bounded coefficients are solved directly. Unlike in \cite{Yo-Z}, the coefficients of equation \eqref{dualeT} are not bounded. Besides equation \eqref{dualeT} we consider the two following equations with values in $\R^{n \times n}$:
\begin{equation}
\left\lbrace
\begin{array}[c]{ll}
 d\Phi_{s}=-H_{s}\Phi_{s}ds+{\displaystyle\sum_{i=1}^{d}}\left( K_{s}^{i}\right) ^{*} \left( K_{s}^{i}\right) ^{*}\Phi_{s}ds-{\displaystyle\sum_{i=1}^{d}}\left( K_{s}^{i}\right) ^{*}\Phi_{s}dW_{s}^{i} & s\in \left[ 0,T\right] \\
\Phi_{0}=I, &
\end{array}
\right.
\label{epPhi}
\end{equation}
and
\begin{equation}
\left\lbrace
\begin{array}[c]{ll}
 d\Psi_{s}=\Psi_{s} H_{s}^{*}ds+{\displaystyle\sum_{i=1}^{d}}\Psi_{s}\left( K_{s}^{i}\right) ^{*}dW_{s}^{i} & s\in \left[ 0,T\right] \\
\Psi_{0}=I. &
\end{array}
\right.  \label{epPsi}
\end{equation}
By applying It\^{o} formula it turns out that $\Phi_{t}\Psi_{t}=I.$
By transposing equation \eqref{epPsi}, we obtain the following equation for $\Psi^{*}$:

\begin{equation}
\left\lbrace
\begin{array}[c]{ll}
 d\Psi_{s}^{*}=H_{s}\Psi_{s}^{*}ds+{\displaystyle\sum_{i=1}^{d}}K_{s}^{i}\Psi_{s}^{*}dW_{s}^{i} & s\in \left[ 0,T\right] \\
\Psi_{0}^{*}=I. &
\end{array}
\right.  \label{epPsi^*}
\end{equation}
By \cite{Gal}, equations \eqref{epPhi}, \eqref{epPsi} and \eqref{epPsi^*} admit a unique solution. Moreover, since $H$ and $K$ are related to the feedback operator, see \eqref{notazionifHK} where $\Lambda$, $H$ and $K$ are defined, it follows that
\begin{equation}
\E \vert \Psi_{t} \vert^{2} \leq C\vert I \vert^{2}, t\in\left[ 0,T \right], \label{stimaPsi}
\end{equation}
where $C$ is a constant that may depend on $T$, see also theorem 2.2 in \cite{Tang}, with $\Psi_{t}^{*}h=\phi_{0,t}h$, $h\in \R^n$.
\noindent We set $\theta:=-\displaystyle\int_{0}^{T}\Psi_{s}P_{s}f_{s}ds$. By boundedness of $P$ and $f$, and by estimate \eqref{stimaPsi} on $\Psi$, it turns out that $\theta \in L^{2}\left( \Omega \right)$. We define
\begin{equation*}
r_{t}=\Phi_{t}[\displaystyle\int_{0}^{t}\Psi_{s}P_{s}f_{s}ds+\Et \theta].
\end{equation*}
Then $r\in L_{loc}^{\infty}\left(\Omega,C([0,T],
\R^{n})\right)$, and following \cite{Yo-Z}, chapter 7, theorem 2.2, we can build a process $g\in L_{loc}^{2}\left(\Omega,C([0,T],
\R^{n})\right)$ such that $(r,g)$ is a solution to equation (\ref{dualeT}), with $(r,g)\in L_{loc}^{\infty}\left(\Omega,C([0,T],
\R^{n})\right) \times L^{2}_{loc}\left(\left[ 0,T\right]
\times\Omega, \R^{n \times d} \right)$.

\finedim

We can prove that the solution $\left( r,g\right) $ to equation \eqref{dualeT} is more regular. To prove this regularity, we need the following \textit{duality relation}.
\begin{remark}\label{lemmadualita}
Let $\xi\in L^2_{\F_T}(\Omega)$ and let $\left( r,g\right) $ be solution to the equation 
\begin{equation*}
\left\lbrace
\begin{array}[c]{ll}
dr_{t}=-H_{t}^{*}r_tdt-P_{t}f_{t}dt-{\displaystyle\sum_{i=1}^{d}\left( K_{t}^{i} \right)^{*}g_{t}^{i}}dt+{\displaystyle\sum_{i=1}^{d}}g_{t}^{i}dW_{t}^{i}, & t \in \left[ 0,T\right]  \\
r_{T}=\xi, &
\end{array}
\right. 
\end{equation*}
and let $X^{t,x,\eta}$ be solution to the equation
\begin{equation}
\left\lbrace
\begin{array}[c]{ll}
 dX_{s}^{t,x,\eta}=H_{s}X_{s}^{t,x,\eta}ds+\displaystyle\sum_{i=1}^{d}K_{s}^{i}X_{s}^{t,x,\eta}dW_{s}^{i}+\eta _{s}ds, & s \in \left[ t,T \right],\\
X_{t}^{t,x ,\eta}=x, & \\
\end{array}
\right.\label{eqdualita}
\end{equation}
where $x\in L^{2}\left( \Omega, \Ft \right)$ and $\eta \in L_{\P}^{2}\left( \Omega \times [0,T], \R^{n} \right)$. Then, by applying the It\^{o} formula, the following \textit{duality relation} holds true:
\begin{equation}
\Et \left\langle \xi, X_{T}^{t,x,\eta}\right\rangle -\left\langle r_{t}, x\right\rangle=-\Et \int_{t}^{T} \left\langle P_{s}f_{s}, X_{s}^{t,x,\eta}\right\rangle ds+\Et \int_{t}^{T} \left\langle \eta_{s}, r_{s}\right\rangle ds.
\label{reldualita}
\end{equation}
\end{remark}

We also need to find a relation between the solution $(r,g)$ of the equation \eqref{dualeT} and the optimal state $\overline{X}$ corresponding to the optimal control $\overline{u}$. This can be achieved, following e.g. \cite{Ben}, by introducing the so called \emph{stochastic Hamiltonian system}
\begin{equation}
\left\{
\begin{array}
[c]{l}
d \overline{X}_{s}=[A_{s}\overline{X}_{s}+B_{s}\overline{u}_{s}]ds+
{\displaystyle\sum_{i=1}^{d}}
[  C_{s}^{i}\overline{X}_{s}+D_{s}^{i}\overline{u} _{s} ] dW_{s}^{i}+f_{s}ds,\\
dy_{s}=-[  A_{s}^{\ast}y_{s}+{\displaystyle\sum_{i=1}^{d}}
\left(  C_{s}^{i}\right)  ^{\ast}z_{s}^{i}+S_{s}\overline{X}_{s}]  ds+
{\displaystyle\sum_{i=1}^{d}}
z_{s}^{i}dW_{s}^{i},\text{ \ \ \ \ \ \ \ \ \ \ \ }t\leq s\leq T,\\
X_{t}=x,\\
y_{T}=\overline{P}_T\overline{X}_{T},
\end{array}
\right.  \label{SHST}%
\end{equation}
where $y,z^{i}\in\mathbb{R}^{n}$, for every $i=1,...,d$. By the so called
stochastic maximum principle, the optimal control for the finite horizon
control problem is given by
\begin{equation}
\overline{u}_{s}=-\left(  B_{s}^{\ast}y_{s}+
{\displaystyle\sum_{i=1}^{d}}
\left(  D_{s}^{i}\right)  ^{\ast}z_{s}^{i}\right)  . \label{u-hamilton}
\end{equation}

By relation (\ref{u-hamilton}), equations (\ref{SHST}) become a fully coupled system of forward backward stochastic differential equations (FBSDE in the following), which admits a unique solution $(\overline{X},y,z) \in  L^{2}_\P\left( \Omega \times [ 0,T], \R^{n} \right)  \times L^{2}_\P\left( \Omega \times [ 0,T], \R^{n}  \right) \times L^{2}_\P\left(\left[ 0,T\right] \times\Omega, \R^{n \times d} \right)$, see Theorem 2.6 in \cite{PeWu}.

\begin{lemma}
Let $(r,g)$ be the unique solution to equation (\ref{dualeT}), and let $(\overline{X},y,z)$ be the unique solution to the FBSDE (\ref{SHST}). Then the following relation holds true for $[0,T]$:
\begin{equation}
y_{t}=P_{t}\overline{X}_{t}+r_{t},\text{ \ \ \ \ \ \ \ \ \ \ \ }t\leq s\leq T.
\label{y-PX-r}
\end{equation}
\label{lemma-y-PX-r}
\end{lemma}
\Dim
We only give a sketch of the proof. For $t=T$ relation (\ref{y-PX-r}) holds true. By applying It\^o formula it turns out that $y_{t}-P_{t}\overline{X}_{t}$ and $r_{t}$ solve the same BSDE, with the same final datum equal to $0$ at the final time $T$. By uniqueness of the solution of this BSDE, the lemma is proved.
\finedim

We can now prove the following regularity result on $(r,g)$.
\begin{proposition}
Let $\left( r,g\right) $ be the solution to equation (\ref{dualeT}). Then $\left( r,g\right) \in L^{2}_\P\left( \Omega,C\left( \left[ 0,T\right], \R^{n} \right) \right) \times L^{2}_\P\left(\left[ 0,T\right] \times\Omega, \R^{n\times d} \right)$. Moreover $r \in L^{\infty}_\P(\Omega \times [0,T])$.
\label{reg_dualeT}
\end{proposition}
\Dim Let $\left( r,g\right) $ be the solution to equation (\ref{dualeT}) built in lemma \ref{esistenzadualeT}.
We note that by theorem 2.6 in \cite{PeWu}, $y \in L^{2}_\P \left( \Omega \times [ 0,T], \R^{n}  \right)$. Moreover, by standard calculations, it is easy to check that $y$ admits a continuous version and $y \in L_\P^{2}\left( \Omega, C([ 0,T], \R^{n} ) \right)$. Moreover, if $f=0$, we get, for every $0 \leq t \leq s \leq T$,
\begin{equation*}
 \Et \vert \overline{X}_{s} \vert ^{2}\leq C \vert x \vert ^{2}.
\end{equation*}
where $\bar X$ is solution to \eqref{SHST}.
This estimate can be easily achieved by applying the Gronwall lemma, and by remembering that from \eqref{controllo_ottimo_Tang_t}, for the optimal control $\overline{u}$ the following holds:
\begin{equation*}
 \Et \int_t^T\vert \overline{u}_{s} \vert ^{2}\leq \langle \overline{P}_t x,x \rangle  \leq C \vert x \vert ^{2}.
\end{equation*}
As a consequence, if $f\neq 0$, for every $0 \leq t \leq s \leq T$,
\begin{equation*}
 \Et \vert \overline{X}_{s} \vert ^{2}\leq C(1+ \vert x \vert ^{2}).
\end{equation*}
Since $P$ is bounded, by lemma \ref{lemma-y-PX-r}, we get that for every $0\leq t \leq s\leq T$
\begin{equation}
 \Et \sup_{t\leq s \leq T} \vert r_{s}\vert ^{2}\leq C,
\label{stima r}
\end{equation}
where $C$ is a constant that can depend on $T$. In particular, for every $0\leq t\leq T$ 
\begin{equation}
 \vert r_{t}\vert ^{2}\leq C.
\label{r-infty}
\end{equation}
Moreover, since $\overline{X}$ is continuous an $P$ admits a continuous version, also $r$ admits a continuous version.
By applying It\^{o} formula we get for $0\leq t\leq w\leq T$,
\begin{align}
\vert r_{t} \vert^{2} & =\vert r_{w} \vert^{2} + 2\int_{t}^{w}\left\langle H_{s} ^{*}r_{s},r_{s}\right\rangle ds+2\int_{t}^{w}\left\langle P_{s}f_{s},r_{s} \right\rangle ds+ 2\int_{t}^{w}\sum_{i=1}^{d}\left\langle \left(  K_{s}^{i}\right) ^{*}g_{s}^{i},r_{s} \right\rangle ds \nonumber \\
& +2 \int_{t}^{w} \sum_{i=1}^{d}\left\langle (g_{s}^{i})^{*} dW_{s}^{i},r_{s} \right\rangle -\int_{t}^{w}\sum_{i=1}^{d}\vert g_{s}^{i} \vert ^{2}ds. 
\label{Ito^2}
\end{align}
By estimate (\ref{r-infty}) and by taking $t=0$, we get by standard calculations
\begin{equation*}
 \E\int_0^{T}\sum_{i=1}^{d}\vert g_{s}^{i} \vert ^{2}ds \leq C,
\end{equation*}
where C is a constant that may depend on $T$. 
\finedim

We are ready to prove the main result of this section
\begin{theorem}
\label{controlloT} Assume $A$, $B$, $C$, $D$ and $f$ satisfy hypothesis \ref{genhypT}. Fix $x\in\mathbb{R}^{n}$, then:
\begin{enumerate}
\item  there exists a unique optimal control $\overline{u}\in L_\P^{2}\left(
\Omega\times\left[  0,T\right]   ,\mathbb{R}^{k}\right)  $ such that for every $0\leq t \leq T$,
\[
J\left(  0,x,\overline{u}\right)  =\inf_{u\in L_\P^{2}\left(
\Omega\times\left[  0,T\right]   ,\mathbb{R}^{k}\right)  }J\left(  0,x,u\right)
\]

\item  If $\overline{X}$ is the mild solution of the state equation
corresponding to $\overline{u}$ (that is the optimal state), then
$\overline{X}$ is the unique mild solution to the \emph{ closed loop }
equation:


\begin{equation}
\left\lbrace
\begin{array}
[c]{ll}
d\overline{X}_{t}=\left[  A_{t}\overline{X}_{t}-B_{t} \left( \Lambda(t, P_t , Q_t)\overline{X}_{t}+( I+
{\displaystyle\sum_{i=1}^{d}}
\left(  D_{t}^{i}\right)  ^{\ast}P_{t}D_{t}^{i})^{-1}(B_{t}^{*}r_{t}+{\displaystyle\sum_{i=1}^{d}}\left( D_{t}^{i}\right) ^{*}g_{t}^{i})\right)\right]  dt + f_t dt & \\
{\displaystyle\sum_{i=1}^{d}}\left[  C_{s}^{i}\overline{X}_{t}-D_{s}^{i}\left(\Lambda(t, P_t , Q_t)\overline{X}_{t}+\left(  I+%

{\displaystyle\sum_{i=1}^{d}}
\left(  D_{t}^{i}\right)  ^{\ast}P_{t}D_{t}^{i}\right)  ^{-1}(B_{t}^{*}r_{t}+{\displaystyle\sum_{i=1}^{d}}\left( D_{t}^{i}\right) ^{*}g_{t}^{i})\right)\right]   dW^i_{t}, &
\\
\overline{X}_{0}=x &
\end{array}
\right.
\label{loopT}%
\end{equation}

\item  The following feedback law holds $\mathbb{P}$-a.s. for almost every
$0 \leq t \leq T$.
\begin{equation}
\overline{u}_{t}=-\left(  I+
\sum_{i=1}^{d}
\left(  D_{t}^{i}\right)  ^{\ast}P_{t}D_{t}^{i}\right)  ^{-1}\left(
P_{t}B_{t}+
\sum_{i=1}^{d}\left(  Q_{t}^{i}D_{t}^{i}+\left(  C_{t}^{i}\right)  ^{\ast}P_{t}D_{t}%
^{i}\right)  \right)  ^{\ast}\overline{X}_{t}+B_{t}^{*}r_{t}+\sum_{i=1}^{d}\left( D_{t}^{i}\right) ^{*}g_{t}^{i}. \label{feedback.T}%
\end{equation}

\item  The optimal cost is given by
\begin{align*}
J(0,x,\overline{u}) & =\langle P_{0}x,x\rangle+2\langle r_0 ,x\rangle-\E \langle P_T \overline{X}_T, \overline{X}_T \rangle +2\E\int_0 ^T \langle r_s ,f_s\rangle ds \nonumber \\
& -\E\int_0 ^T \vert (  I+
{\displaystyle\sum_{i=1}^{d}}
\left(  D_{t}^{i}\right)  ^{\ast}P_{t}D_{t}^{i})  ^{-1}( B_{t}^{*}r_{t}+{\displaystyle\sum_{i=1}^{d}}\left( D_{t}^{i}\right) ^{*}g_{t}^{i})\vert ^2 ds.
\end{align*}
\end{enumerate}
\end{theorem}

\noindent\hbox{{\bf Proof.}$\;\; $} By computing $d\langle P_t X_t, X_t \rangle +2\langle r_t , X_t\rangle$, we get the so called fundamental relation
\begin{align*}
& \Et \int_{t}^{T}[\left\langle S_{s}X_{s}^{0,x,u},X_{s}^{0,x,u}\right\rangle +|u_{s}|^{2}]ds \nonumber \\
& =\langle P_{t}x,x\rangle+2\langle r_t ,x\rangle-\Et \langle P_T X_T, X_T \rangle +2\Et\int_t ^T \langle r_s ,f_s\rangle ds \nonumber \\
&+\Et \!\!\int_t ^T \!\!\vert \!\!\left(  I+ \sum_{i=1}^{d} \left(
D_{s}^{i}\right)  ^{\ast}P_{s}D_{s}^{i}\right)
^{-1}\!\!\!\!\!\left( P_{s}B_{s}+
\sum_{i=1}^{d}\left(  Q_{s}^{i}D_{s}^{i}+\left(  C_{s}^{i}\right)  ^{\ast}P_{s}D_{s}%
^{i}\right)  \right)  ^{\ast}\!\!\!X_{s}+B_{s}^{*}r_{s}+\sum_{i=1}^{d}\left( D_{s}^{i}\right)\vert^2 ds \nonumber \\
&-\Et\int_t ^T \vert (  I+\sum_{i=1}^{d}\left(  D_{s}^{i}\right)  ^{\ast}P_{s}D_{s}^{i})  ^{-1}(B_{s}^{*}r_{s}+{\displaystyle\sum_{i=1}^{d}}\left( D_{s}^{i}\right) ^{*}g_{s}^{i})\vert ^2 ds.
\end{align*}
The theorem now easily follows.
{\hfill\hbox{\enspace${ \square}$}}

\section{Preliminary results for the infinite horizon case}

The next step is to study the optimal control problem in the infinite horizon case and with $f\neq 0 $. To this aim we have to study solvability and regularity of the solution of a BSRDE with infinite horizon, in particular we study $P$. At first we consider the case when $f=0$. Namely, in this section we consider the following stochastic differential equation where $X^{t,x,u}$ represents the state:
\begin{equation}
\left\{
\begin{array}
[c]{ll}
dX_{s}^{t,x,u}=(A_{s}X_{s}^{t,x,u}+B_{s}u_{s})ds+
{\displaystyle\sum_{i=1}^{d}}
\left(  C_{s}^{i}X_{s}^{t,x,u}+D_{s}^{i}u_{s}\right)  dW_{s}^{i} & s\geq t\\
X_{t}^{t,x,u}=x &
\end{array}
\right.  \label{stato.stato}%
\end{equation}
As a by product of the preliminaries studies, we are able to solve the following stochastic optimal control problem: minimize with respect to every admissible control $u$ the cost functional,
\begin{equation}
J_{\infty}(0,x,u)=\mathbb{E}\int_{0}^{+\infty}[\left\langle S_{s}X_{s}%
^{0,x,u},X_{s}^{0,x,u}\right\rangle +|u_{s}|^{2}]ds.
\end{equation}
We define the set of admissible control
\begin{equation}\label{U_ad}
 \mathcal U =\left\lbrace u\in L^2_\P(\Omega\times[0,+\infty),\R^k):\mathbb{E}\int_{0}^{+\infty}\left\langle S_{s}X_{s}%
^{0,x,u},X_{s}^{0,x,u}\right\rangle +|u_{s}|^{2}ds< +\infty\right\rbrace .
\end{equation}
We also introduce the following random variables, for $t\in\lbrack0,+\infty]$:%

\[
J_{\infty}(t,x,u)=\mathbb{E}^{\mathcal{F}_{t}}\int_{t}^{+\infty}[\left\langle
S_{s}X_{s}^{t,x,u},X_{s}^{t,x,u}\right\rangle +|u_{s}|^{2}]ds
\]

We will work under the following general assumptions on $A$, $B$, $C$ and
$D$\ that will hold from now on:

\begin{hypothesis}
\label{genhyp} $\ $

\begin{enumerate}
\item [ A1)]$A:\left[  0,+\infty\right)  \times\Omega\rightarrow$
$\mathbb{R}^{n\times n}$, $B:\left[  0,+\infty\right)
\times\Omega \rightarrow$ $\mathbb{R}^{n\times k}$, $C^{i}:\left[
0,+\infty\right) \times\Omega\rightarrow$ $\mathbb{R}^{n\times
n}$, $i=1,...,d$ and $D^{i}:\left[  0,+\infty\right)
\times\Omega\rightarrow$ $\mathbb{R}^{n\times k}$, $i=1,...,d$,
are uniformly bounded process adapted to the filtration $\left\{
\mathcal{F}_{t}\right\}  _{t\geq0}$.

\item[ A2)] $S:\left[  0,+\infty\right)  \times\Omega\rightarrow
\mathbb{R}^{n\times n}$ is uniformly bounded and adapted to the
filtration $\left\{  \mathcal{F}_{t}\right\}  _{t\geq0}$ and it is
almost surely and almost everywhere symmetric and nonnegative.
\end{enumerate}
\end{hypothesis}

In order to study this control problem in infinite horizon, we consider the following backward stochastic
Riccati equation on $[0,+\infty)$:%
\begin{align}
&dP_{t}    =-\left[  A_{t}^{\ast}P_{t}+P_{t}A_{t}+S_{t}+
{\displaystyle\sum_{i=1}^{d}}
\left(  \left(  C_{t}^{i}\right)  ^{\ast}P_{t}C_{t}^{i}+\left(  C_{t}
^{i}\right)  ^{\ast}Q_{t}+Q_{t}C_{t}^{i}\right)  \right]  dt+
{\displaystyle\sum_{i=1}^{d}}
Q_{t}^{i}dW_{t}^{i}+\label{RiccatiAlg}\\
&  \left[  P_tB_t+ {\displaystyle\sum_{i=1}^{d}} \left(  \left(
C^{i}_t\right)  ^{\ast}P_t D^{i}_t +Q^{i}D_t ^{i}\right)  \right]
\left[  I+ {\displaystyle\sum_{i=1}^{d}} \left(  D^{i}_t\right)
^{\ast}P_tD^{i}_t\right]  ^{-1}\!\!\!\left[  P_tB_t+
{\displaystyle\sum_{i=1}^{d}} \left(  \left(  C^{i}_t\right)
^{\ast}P_t D_t ^{i}+Q^{i}_t D^{i}_t\right) \right]
^{\ast}\!\!\!\!\, dt,\nonumber
\end{align}
where we stress that the final condition has disappeared but we ask that the
solution can be extended to the whole positive real half-axis.

\begin{definition}\label{solriccati_inf}
We say that a pair of processes $\left(  P,Q\right)  $ is a solution to
equation (\ref{RiccatiAlg}) if for every $T>0$ $\left(  P,Q\right)  $ is a
solution to equation (\ref{Riccati0T})\ in the interval time $\left[
0,T\right]  $, with  with final datum the process $P$ evalueted at time $T$.
\noindent In particular a solution $(\overline{P},\overline{Q})$ is called minimal if 
whenever another couple $(P,Q)$ is a solution to the Riccati equation then
$P -\overline{P}$ is a non-negative matrix, see also Corollary 3.3 in
\cite{GTinf}.
\end{definition}

\begin{definition}
\label{defstab}We say that $(A,B,C,D)$ is stabilizable relatively to the
observations $\sqrt{S}$ (or $\sqrt{S}$-stabilizable) if there exists a control
$u\in L_{\mathcal{P}}^{2}([0,+\infty)\times\Omega;U)$ such that for all
$t\geq0$ and all $x\in\mathbb{R}^{n}$%
\begin{equation}
\mathbb{E}^{\mathcal{F}_{t}}\int_{t}^{+\infty}[\left\langle S_{s}X_{s}%
^{t,x,u},X_{s}^{t,x,u}\right\rangle +|u_{s}|^{2}]ds<M_{t,x}.\label{condstabS}%
\end{equation}
for some positive constant $M_{t,x}$.
\end{definition}
\noindent This kind of stabilizability condition has been
introduced in \cite{GTinf}.

In the following, we consider BSRDEs on the time interval $\left[
0,N\right] $, with final condition $P_{N}=0$. For each integer
$N>0$, let $\left(
P^{N},Q^{N}\right)  $ be the solution of the Riccati equation%
\begin{equation}\label{RiccatiN}
\left\{
\begin{array}
[c]{l}
-dP_{t}^{N}=G\left(  A_{t},B_{t},C_{t},D_{t};S_{t};P_{t}^{N},Q_{t}^{N}\right)
dt+
{\displaystyle\sum_{i=1}^{d}}
Q_{t}^{N,i}dW_{t}^{i}\\
P_{N}^{N}=0.
\end{array}
\right.
\end{equation}
$P^{N}$ can be defined in the whole $[0,+\infty)$ setting $P_{t}^{N}=0$ for
all $t>N$. We prove the following lemma.

\begin{lemma}\label{convergenzaN}
Assume hypothesis \ref{genhyp} and that $(A,B,C,D)$ is stabilizable relatively to the
observations $\sqrt{S}$. There exists a random matrix $\overline{P}$\ almost
surely positive and symmetric such that $\ \mathbb{P}\left\{  \lim
_{N\rightarrow+\infty}P^{N}(t)x=\overline{P}(t)x,\ \forall x\in\mathbb{R}%
^{n}\right\}  =1$.
\end{lemma}

\noindent\hbox{{\bf Proof.}$\;\; $}The proof essentially follows
the first part of the proof of proposition 3.2 in \cite{GTinf}.
\finedim

\begin{remark}
It is clear from the above proof of proposition 3.2 in \cite{GTinf} that condition \eqref{condstabS}
is equi\-va\-lent to the following one:
\begin{equation}
\mathbb{E}^{\mathcal{F}_{t}}\int_{t}^{+\infty}[\left\langle S_{s}X_{s}%
^{t,x,u},X_{s}^{t,x,u}\right\rangle +|u_{s}|^{2}]ds<M |x|^2.\label{condstabSnorm}%
\end{equation}
where the constant $M$ may depend on $t$.
\end{remark}
Next we want to prove that $\overline{P}$ built in the previous lemma is the solution to the BSRDE (\ref{RiccatiAlg}). This is achieved through the control meaning of the solution of the Riccati equation. Indeed also the martigale term $Q^N$ appears in the fundamental relation is no more possible to proceed as in \cite{GTinf}.
For $T>0$ fixed and for each $N>T$, we consider the following finite horizon
stochastic optimal control problem: minimize the cost, over all admissible controls,
\begin{equation*}
J(0,x,u)=\mathbb{E}\langle
P^{N}_TX_{T}^{0,x,u},X_{T}^{0,x,u}\rangle
+\mathbb{E}\int_{0}^{T}[\left\langle S_{s}X_{s}^{0,x,u},X_{s}^{0,x,u}%
\right\rangle +|u_{s}|^{2}]ds, 
\end{equation*}
where $X^{0,x,u}$\ is solution to equation (\ref{stato.stato}). Let $u^{N}$ be the optimal
control, and $X^{N}$ the corresponding optimal state. Let $\widetilde{u}$ be the
optimal control, and $\widetilde{X}$ the corresponding optimal state for the following finite horizon optimal control problem: minimize the cost, over all admissible controls,
\begin{equation*}
J(0,x,u)=\mathbb{E}\langle\overline{P}_TX_{T}^{0,x,u},X_{T}^{0,x,u}\rangle
+\mathbb{E}\int_{0}^{T}[\left\langle S_{s}X_{s}^{0,x,u},X_{s}^{0,x,u}%
\right\rangle +|u_{s}|^{2}]ds. 
\end{equation*}
Let us consider the so called \emph{stochastic Hamiltonian system}
\begin{equation}
\left\{
\begin{array}
[c]{l}
dX_{s}=[A_{s}X_{s}-B_{s}(  B_{s}^{\ast}y_{s}+
{\displaystyle\sum_{i=1}^{d}}
\left(  D_{s}^{i}\right)  ^{\ast}z_{s}^{i})]  ds+
{\displaystyle\sum_{i=1}^{d}}
[  C_{s}^{i}X_{s}+D_{s}^{i}( B_{s}^{\ast}y_{s}+
{\displaystyle\sum_{k=1}^{d}}
 (D_{s}^{k})^{\ast}z_{s}^{k} ) ] dW_{s}^{i},\\
dy_{s}=-[  A_{s}^{\ast}y_{s}+{\displaystyle\sum_{i=1}^{d}}
\left(  C_{s}^{i}\right)  ^{\ast}z_{s}^{i}+S_{s}X_{s}]  ds+
{\displaystyle\sum_{i=1}^{d}}
z_{s}^{i}dW_{s}^{i},\text{ \ \ \ \ \ \ \ \ \ \ \ }t\leq s\leq T,\\
X_{t}=x,\\
y_{T}=\overline{P}_TX_{T},
\end{array}
\right.  \label{SHS}%
\end{equation}
where $y,z^{i}\in\mathbb{R}^{n}$, for every $i=1,...,d$. By the so called
stochastic maximum principle, the optimal control of the finite horizon
control problem is given by
\[
u_{s}=-\left(  B_{s}^{\ast}y_{s}+
{\displaystyle\sum_{i=1}^{d}}
\left(  D_{s}^{i}\right)  ^{\ast}z_{s}^{i}\right)  .
\]

Let us consider the stochastic Hamiltonian systems relative to the optimal
control $u^{N}$ and to the optimal control $\widetilde{u}$, and let us denote
by $\left(  X^{N},y^{N},z^{N}\right)  $ and by $\left(  \widetilde
{X},\widetilde{y},\widetilde{z}\right)  $ the solutions of the corresponding
stochastic Hamiltonian systems.

\begin{lemma}
\label{lemmaSHS}$\mathbb{E}^{\mathcal{F}_{t}}\displaystyle\int_{t}^{T}[\left|  \sqrt{S_{s}%
}(\widetilde{X}_{s}-X_{s}^{N})\right| ^2 +|B_{s}^{\ast}(\widetilde{y}_{s}%
-y_{s}^{N})+
{\displaystyle\sum_{i=1}^{d}}
\left(  D_{s}^{i}\right)  ^{\ast}(\widetilde{z}_{s}^{i}-z_{s}^{N,i}%
)|^{2}]ds\rightarrow0$ as $N\rightarrow\infty$.
\end{lemma}

\noindent\hbox{{\bf Proof.}$\;\; $} The proof is based on the
application of It\^ o formula to $\langle \tilde{y}_t -
y_{t}^{N},\widetilde{X}_{t}-X_{t}^{N} \rangle$.
\begin{align*}
&  \mathbb{E}^{\mathcal{F}_{t}}\langle \widetilde{y}_{T}-y_{T}
^{N},\widetilde{X}_{T}-X_{T}^{N}\rangle =
 \mathbb{E}^{\mathcal{F}_{t}}\int_{t}^{T}d\langle \widetilde{y}
_{s}-y_{s}^{N},\widetilde{X}_{s}-X_{s}^{N}\rangle \\
 & =-\mathbb{E}^{\mathcal{F}_{t}}
\int_{t}^{T}[  \left|  B_{s}^{\ast}\left(  \widetilde{y}_{s}-y_{s}%
^{N}\right)  \right|^{2}+\left| {\displaystyle\sum_{i=1}^{d}}
\left(  D_{s}^{i}\right)^{\ast}\left(  \widetilde{z}_{s}^{i}-z_{s}%
^{N,i}\right)  \right|  ^{2}+ 2 {\displaystyle\sum_{i=1}^{d}}
\langle B_{s}^{\ast}\left( \widetilde{y}_{s}-y_{s}^{N}\right)
,\left( D_{s}^{i}\right) ^{\ast}\left(
\widetilde{z}_{s}^{i}-z_{s}^{N,i}\right) \rangle ] ds
\\
& -\mathbb{E}^{\mathcal{F}_{t}} {\displaystyle\sum_{i=1}^{d}}
\int_{t}^{T}\left| \sqrt{S_{s}}\left(
\widetilde{X}_{s}-X_{s}^{N}\right) \right| ^{2}ds\\ &  =
-\mathbb{E}^{\mathcal{F}_{t}} {\displaystyle\sum_{i=1}^{d}}
\int_{t}^{T}\left|
B_{s}^{\ast}\left(\widetilde{y}_{s}-y_{s}^{N}\right)  +
{\displaystyle\sum_{i=1}^{d}} \left(D_{s}^{i}\right) ^{\ast}\left(
\widetilde{z}_{s}^{i}-z_{s} ^{N,i}\right) \right|^{2}ds
-\mathbb{E}^{\mathcal{F}_{t}} {\displaystyle\sum_{i=1}^{d}}
\int_{t}^{T}\left|\sqrt{S_{s}}\left(
\widetilde{X}_{s}-X_{s}^{N}\right) \right|^{2}ds.
\end{align*}
Since $\widetilde{y}_{T}=\overline{P}_{T}\widetilde{X}_{T}$ and $y_{T}
^{N}=P_{T}^{N}X_{T}^{N}$, we finally get
\begin{align*}
&  \mathbb{E}^{\mathcal{F}_{t}}\left\langle \overline{P}_{T}\widetilde{X}
_{T}-P_{T}^{N}X_{T}^{N},\widetilde{X}_{T}-X_{T}^{N}\right\rangle \\
&  =-\mathbb{E}^{\mathcal{F}_{t}}
{\displaystyle\sum_{i=1}^{d}}
\int_{t}^{T}\left|  B_{s}^{\ast}\left(  \widetilde{y}_{s}-y_{s}^{N}\right)  +
{\displaystyle\sum_{i=1}^{d}}
\left(  D_{s}^{i}\right)  ^{\ast}\left(  \widetilde{z}_{s}^{i}-z_{s}
^{N,i}\right)  \right|  ^{2}ds-\mathbb{E}^{\mathcal{F}_{t}}
{\displaystyle\sum_{i=1}^{d}}
\int_{t}^{T}\left|  \sqrt{S_{s}}\left(  \widetilde{X}_{s}-X_{s}^{N}\right)
\right|  ^{2}ds.
\end{align*}
By adding and subtracting $\mathbb{E}^{\mathcal{F}_{t}}\left\langle P_{T}
^{N}\widetilde{X}_{T},\widetilde{X}_{T}-X_{T}^{N}\right\rangle ,$
\begin{align*}
&  \mathbb{E}^{\mathcal{F}_{t}}\left\langle P_{T}^{N}\left(  \widetilde{X}
_{T}-X_{T}^{N}\right)  ,\widetilde{X}_{T}-X_{T}^{N}\right\rangle
+\mathbb{E}^{\mathcal{F}_{t}}\left\langle \left(  \overline{P}_{T}-P_{T}
^{N}\right)  \widetilde{X}_{T},\widetilde{X}_{T}-X_{T}^{N}\right\rangle \\
&  =-\mathbb{E}^{\mathcal{F}_{t}}
\int_{t}^{T}\left|  B_{s}^{\ast}\left(  \widetilde{y}_{s}-y_{s}^{N}\right)  +
\left(  D_{s}^{i}\right)  ^{\ast}\left(  \widetilde{z}_{s}^{i}-z_{s}
^{N,i}\right)  \right|  ^{2}ds-\mathbb{E}^{\mathcal{F}_{t}}
{\displaystyle\sum_{i=1}^{d}}
\int_{t}^{T}\left|  \sqrt{S_{s}}\left(  \widetilde{X}_{s}-X_{s}^{N}\right)
\right|  ^{2}ds.
\end{align*}
Since $\left\langle P_{T}^{N}\left(  \widetilde{X}_{T}-X_{T}^{N}\right)
,\widetilde{X}_{T}-X_{T}^{N}\right\rangle \geq0$, and by definition
$\mathbb{E}^{\mathcal{F}_{t}}\left\langle \left(  \widetilde{P}_{T}-P_{T}
^{N}\right)  \widetilde{X}_{T},\widetilde{X}_{T}-X_{T}^{N}\right\rangle
\rightarrow0$ for $N$ sufficiently large, also
\[
\mathbb{E}^{\mathcal{F}_{t}}
\int_{t}^{T}\left|  B_{s}^{\ast}\left(  \widetilde{y}_{s}-y_{s}^{N}\right)  +
{\displaystyle\sum_{i=1}^{d}}
\left(  D_{s}^{i}\right)  ^{\ast}\left(  \widetilde{z}_{s}^{i}-z_{s}
^{N,i}\right)  \right|  ^{2}ds+\mathbb{E}^{\mathcal{F}_{t}}
\int_{t}^{T}\left|  \sqrt{S_{s}}\left(  \widetilde{X}_{s}-X_{s}^{N}\right)
\right|  ^{2}ds\rightarrow0
\]
as $N\rightarrow\infty$. In particular this means that $\mathbb{E}
^{\mathcal{F}_{t}}\displaystyle\int_{t}^{T}\left|  \widetilde{u}_{s}-u_{s}^{N}\right|
^{2}ds\rightarrow 0$ as $N\rightarrow\infty$.{\hfill\hbox{\enspace${ \square}$}}

As a consequence of the previous results we deduce the following:
\begin{corollary}\label{CorMin}
Assume hypothesis \ref{genhyp} and that $(A,B,C,D)$ is stabilizable relatively to $\sqrt{S}$. The process $\overline{P}$ is the minimal solution of the Riccati equation in the sense of definition \ref{solriccati_inf}.
\end{corollary}
\Dim
Fix $T<N$ and on $[0,T]$ consider the Riccati equation
\begin{equation}
\left\{
\begin{array}
[c]{l}
-dP^N_{t}=G\left(  A_{t},B_{t},C_{t},D_{t};S_{t};P^N_{t},Q^N_{t}\right)  dt+
{\displaystyle\sum_{i=1}^{d}}
Q_{t}^{i,N}dW_{t}^{i}\\
P^N_T=P^N(T).
\end{array}
\right.
\end{equation}
Then \begin{equation}
\left\langle P_{t}^{N}x,x\right\rangle    =\mathbb{E}^{\mathcal{F}_{t}
}\left\langle P_{T}^{N}X_{T}^{N},X_{T}^{N}\right\rangle +\mathbb{E}
^{\mathcal{F}_{t}}\int_{t}^{T}[\left\langle S_{r}X_{r}^{N},X_{r}
^{N}\right\rangle +|u^N_{r}|^{2}]dr\label{P_Ncor}
\end{equation}
By lemma \ref{lemmaSHS} we deduce that
\[\mathbb{E}^{\mathcal{F}_{t}}\left\langle P_{T}^{N}\left(  \widetilde{X}
_{T}-X_{T}^{N}\right)  ,\widetilde{X}_{T}-X_{T}^{N}\right\rangle \to 0 \qquad \text{ as } N \to +\infty.\]
So $\mathbb{E}^{\mathcal{F}_{t}
}\left\langle P_{T}^{N}X_{T}^{N},X_{T}^{N}\right\rangle \to \mathbb{E}^{\mathcal{F}_{t}
}\left\langle\overline{P}_{T}\widetilde{X}_{T},\widetilde{X}_{T}\right\rangle $ as $N\rightarrow +\infty$, since
\begin{equation*}
\mathbb{E}^{\mathcal{F}_{t}
}\left\langle P_{T}^{N}X_{T}^{N},X_{T}^{N}\right\rangle=\mathbb{E}^{\mathcal{F}_{t}}\left\langle P_{T}^{N}\left(  \widetilde{X}_{T}-X_{T}^{N}\right)  ,\widetilde{X}_{T}-X_{T}^{N}\right\rangle+2\langle P_{T}^{N}\widetilde{X}_{T},X_{T}^{N} \rangle-\langle P_{T}^{N}\widetilde{X}_{T},\widetilde{X}_{T} \rangle.
\end{equation*}
By the construction of $\overline{P}$ and by lemma \ref{lemmaSHS}, we have that, letting $N\rightarrow +\infty$ in \eqref{P_Ncor} \begin{equation*}
\left\langle \overline{P}_{t}x,x\right\rangle    =\mathbb{E}^{\mathcal{F}_{t}
}\left\langle\overline{P}_{T}\widetilde{X}_{T},\widetilde{X}_{T}\right\rangle +\mathbb{E}
^{\mathcal{F}_{t}}\int_{t}^{T}[\left\langle S_{r}\widetilde{X}_{r},\widetilde{X}_{r}\right\rangle +|\widetilde{u}_{r}|^{2}]dr
\end{equation*}
So $(\overline{P},\overline{Q})$ is the minimal solution of the Riccati equation in the sense of definition \ref{solriccati_inf}.
\finedim


By the previous calculations, we can now solve the optimal control problem with infinite horizon, when $f=0$.
\begin{theorem}
\label{controllo}If $A1)-A2)$ hold true and if $(A,B,C,D)$ is stabilizable
relatively to $S$, given $x\in\mathbb{R}^{n}$, then:
\begin{enumerate}
\item  there exists a unique optimal control $\overline{u}\in L^{2}_\P\left(
\Omega\times\left[  0,+\infty\right)  ,\mathbb{R}^{k}\right)  $ such that
\[
J_{\infty}\left(  0,x,\overline{u}\right)  =\inf_{u\in L^{2}_\P\left(
\Omega\times\left[  0,+\infty\right)  ,\mathbb{R}^{k}\right)  }J_{\infty
}\left(  0,x,u\right)
\]

\item The process $\overline{P}$ defined in lemma 3.4 is the minimal solution of the Riccati equation.

\item  If $\overline{X}$ is the mild solution of the state equation
corresponding to $\overline{u}$ (that is the optimal state), then
$\overline{X}$ is the unique mild solution to the \emph{ closed loop }
equation:
\begin{equation}
\left\{
\begin{array}
[c]{ll}
d\overline{X}_{t}=\left[  A\overline{X}_{t}-B_{t}\left(  I+%

{\displaystyle\sum_{i=1}^{d}}
\left(  D_{t}^{i}\right)  ^{\ast}\overline{P}_{t}D_{t}^{i}\right)  ^{-1}\left(
\overline{P}_{t}B_{t}+
{\displaystyle\sum_{i=1}^{d}}
\left(  \overline{Q}_{t}^{i}D_{t}^{i}+\left(  C_{t}^{i}\right)  ^{\ast}\overline{P}_{t}D_{t}%
^{i}\right)  \right)  ^{\ast}\overline{X}_{t}\right]  dt +& \\
{\displaystyle\sum_{i=1}^{d}}\left[  C_{t}^{i}\overline{X}_{t}-D_{t}^{i}\left(  I+%

{\displaystyle\sum_{i=1}^{d}}
\left(  D_{t}^{i}\right)  ^{\ast}\overline{P}_{t}D_{t}^{i}\right)  ^{-1}\left(
\overline{P}_{t}B_{t}+
{\displaystyle\sum_{i=1}^{d}}

\left(  \overline{Q}_{t}^{i}D_{t}^{i}+\left(  C_{t}^{i}\right)  ^{\ast}\overline{P}_{t}D_{t}%
^{i}\right)  \right)  ^{\ast}\overline{X}_{t}\right]   dW_{t}, &
\\
\overline{X}_{0}=x. &
\end{array}
\right.  \label{loopinf}%
\end{equation}

\item  The following feedback law holds $\mathbb{P}$-a.s. for almost every
$t$:
\begin{equation}
\overline{u}_{t}=-\left(  I+
{\displaystyle\sum_{i=1}^{d}}
\left(  D_{t}^{i}\right)  ^{\ast}\overline{P}_{t}D_{t}^{i}\right)  ^{-1}\left(
\overline{P}_{t}B_{t}+
{\displaystyle\sum_{i=1}^{d}}\left( \overline{Q}_{t}^{i}D_{t}^{i}+\left(  C_{t}^{i}\right)  ^{\ast}\overline{P}_{t}D_{t}%
^{i}\right)  \right)  ^{\ast}\overline{X}_{t}. \label{feedback.ban.inf}%
\end{equation}

\item  The optimal cost is given by $J_{\infty}(0,x,\overline{u})=\langle
\overline{P}_{0}x,x\rangle$.
\end{enumerate}
\end{theorem}
The proof of this theorem is similar, and more immediate, to the proof of theorem \ref{teocontrollo}, which is given in detail in section 5. In particular we deduce that
\begin{equation}
\langle\overline{P}_{t}x,x\rangle=\mathbb{E}^{\mathcal{F}_{t}}\int_{t}^{\infty
}[\left\langle S_{r}\widetilde{X}_{r}^{t,x,u},\widetilde{X}_{r}^{t,x,u}%
\right\rangle +|\widetilde{u}_{r}|^{2}]dr.\label{u_lim}
\end{equation}

\section{The infinite horizon dual equation}
We first introduce some definitions. We say that a solution $P$ of equation \ref{RiccatiAlg}
is  bounded, if there exists a constant $M >0$  such that for
every $t \geq 0$
\begin{equation*}
|P_t| \leq M \qquad\qquad \mathbb{P}-a.s.
\end{equation*}
Whenever the constant $M_{t,x}$ that appears in definition
\ref{defstab} can be chosen independently of $t$, then the minimal
solution $\overline{P}$ is automatically bounded.
\begin{definition}\label{def.stab.identita}
Let $P$ be a solution to \ref{RiccatiAlg}. We say that $P$
stabilizes $(A,B,C,D)$ relatively to the identity $I$ if for every
$t > 0$ and $x \in \R^n$ there exists a positive constant $M$,
independent of $t$, such that
\begin{equation}\label{condstabiden}
\Et\int_t^{+\infty}|X^{t,x}_r|^2, dr \leq M \quad \quad
\mathbb{P}-a.s.,\end{equation} where $X^{t,x}$ is, for $s\geq t$, the
solution to:
\begin{equation}
\left\{
\begin{array}
[c]{ll}
dX_s^{t,x}=H_s



X_s^{t,x}  ds +

{\displaystyle\sum_{i=1}^{d}}

K^i_sX_s^{t,x}   dW^i_{s}, &
\\
X_t^{t,x}=x &
\end{array}
\right.\label{closedloop_tx}
\end{equation}
\end{definition}
From now on we assume that the process $\overline{P}$ is bounded and stabilizes $(A,B,C,D)$ with respect to the identity $I$.
\begin{remark}
 It is possible to verify in some concrete situations that $(A,B,C,D)$ is stabilizable relatively to
  the observations $\sqrt{S}$ and that $\overline{P}$ stabilizes $(A,B,C,D)$ relatively to the identity $I$.
   Here we present the case when, for some $\alpha >0$, $A$ and $C$ satisfy
\begin{equation}
 \langle A_tx,x\rangle+\frac{1}{2}\langle C_tx,C_tx\rangle\leq -\alpha \abs{x}^2, \label{cond_A-C}
\end{equation}
for every $t\geq 0$ and $x\in \R^n$, then $(A,B,C,D)$ is stabilizable relatively to the observations $\sqrt{S}$
 uniformly in time. Indeed, by taking the control $u=0$, applying the It\^{o} formula to the state equation we get, for $0\leq t\leq s$,
\begin{align*}
\Et \abs{X_s^{t,x,0}}^2 &\leq \abs{x}^2+2\Et \int_t^s \langle A_rX_r^{t,x,0},X_r^{t,x,0}\rangle dr +\Et \int_t^s \langle C_rX_r^{t,x,0},C_rX_r^{t,x,0}\rangle dr\\
&\leq \abs{x}^2-2\alpha\Et \int_t^s \abs{X_r}^2dr.
\end{align*}
By the Gronwall lemma we get
\begin{equation*}
\Et \abs{X_s^{t,x,0}}^2 \leq \abs{x}^2 e^{-2\alpha (s-t)}.
\end{equation*}
So for every $0\leq t\leq T$
\begin{equation*}
\Et \int_t^{+\infty}\abs{\sqrt{S}X_r^{t,x,0}}^2 dr \leq M_x,
\end{equation*}
where $M_x$ is a constant dependent on the initial condition $x$, but independent on the initial time $t$. So,
 according to definition \ref{defstab}, $(A,B,C,D)$ is stabilizable relatively to the observations $\sqrt{S}$,
  uniformly in time. Moreover, assuming that $S\geq \epsilon I$, for some $\epsilon >0$, by \eqref{cond_A-C},
   we also get that $\overline{P}$ stabilizes $(A,B,C,D)$ relatively to the identity $I$. Indeed, by the previous calculations, denoting by $\overline X$ and $\overline u$ respectively the optimal state and the optimal control for the infinite horizon control problem with $f=0$, it follows that
\begin{equation*}
\Et \int_t^{+\infty}[\abs{\overline X_r}^2+\abs{\overline u}^2] dr \leq\Et \int_t^{+\infty}\abs{X_r^{t,x,0}}^2 dr \leq M_x,
\end{equation*}
which immediately implies \eqref{condstabiden}.
\end{remark}

The Datko theorem holds also in this case. 
\begin{theorem}\label{remarkDatko}
If $\overline{P}$ stabilizes $(A,B,C,D)$ relatively to the identity
and it is uniformly bounded in time, then for the solution of equation (\ref{closedloop_tx}) there exists 
two constants $a>0$, $C>0$ such that 

\begin{equation}\label{stimadatkoinf}
\Et \vert X^{t,x}_s\vert^2 \leq C e^{-a(s-t)}\vert x \vert ^2\qquad \mathbb{P}-\text{a.s.}
\end{equation}
\end{theorem}

\begin{remark}Note that if $X_s^{t,x}$ is solution of equation
 \begin{equation*}
\left\{
\begin{array}
[c]{ll}
dX_s^{t,x,\eta}=H_sX_s^{t,x,\eta}+

{\displaystyle\sum_{i=1}^{d}}

K^i_sX_s^{t,x,\eta}   dW^i_{s}+\eta_sds, &
\\
X_t^{t,x,\eta}=x, &
\end{array}
\right.
\end{equation*} then by the previous theorem
\[
 \int_0^{+\infty}\E\vert X_s^{t,x,\eta} \vert^2 ds \leq C \left(\int_0^{+\infty} \E \vert\eta_s\vert^2
+\vert x\vert^2 \right)
\label{stimanonlin}
\]

\end{remark}

In order to study the optimal control problem with infinite horizon and with $f\neq 0$, we need to study the BSDE on $[0,\infty)$,
\begin{equation}
dr_{t}=-H_{t}^{*}r_t dt-P_{t}f_{t}dt-{\displaystyle\sum_{i=1}^{d}\left( K_{t}^{i} \right)^{*}g_{t}^{i}}dt+{\displaystyle\sum_{i=1}^{d}}g_{t}^{i}dW_{t}^{i},\text{
\ \ \ \ \ } t \geq 0  \\
\label{duale}
\end{equation}
where the final condition has disappeared but we ask that the solution can be extended to the whole positive real axis.
We make the following assumption on $f$:
\begin{hypothesis}
$f$ is a process in $L^2_\P (\Omega \times [0,+\infty), \R^n) \cap L^{\infty}_\P(\Omega \times [0,+\infty), \R^{n }) $.
\label{hypf}
\end{hypothesis}
\begin{proposition}
Let hypotheses \ref{genhyp} and \ref{hypf} hold true and assume that $\overline{P}$ is bounded and stabilizes $(A,B,C,D)$ with respect to the identity $I$. Then equation (\ref{duale}) admits a solution $(\bar r, \bar g) \in  L^2_\P (\Omega \times [0,+\infty), \R^n)\times L^2_\P(\Omega \times [0,T], \R^{n \times d}) $, for every $T>0$.
\label{reg-r-infty}
\end{proposition}
\Dim
For integer $N>0$, we consider the BSDEs
\begin{equation}
\left\lbrace
\begin{array}[c]{ll}
dr^N _{t}=-H_{t}^{*}r_t ^N dt-P_{t}f_{t}dt-{\displaystyle\sum_{i=1}^{d}\left( K_{t}^{i} \right)^{*}g_{t}^{i,N}}dt+{\displaystyle\sum_{i=1}^{d}}g_{t}^{i,N}dW_{t}^{i}, & t \in \left[ 0,T\right]  \\
r_N ^N=0. &
\end{array}
\right. \label{dualeN}
\end{equation}
By proposition \ref{reg_dualeT}, we know that equation
(\ref{dualeN}) admits a unique solution $\left( r^N,g^N\right)$
that belongs to $L^{2}_\P\left( \Omega,C\left( \left[ 0,N\right],
\R^{n} \right) \right) \times L^{2}_\P\left(\left[ 0,N\right]
\times\Omega, \R^{n\times d} \right)$, for every $N\in \nat$. The
aim is to write a \textit{duality} relation, see remark
\ref{lemmadualita}, between $r^N$ and the process $X^N$, solution
of the following equation
\begin{equation*}
\left\lbrace
\begin{array}[c]{ll}
dX^N _{s}=H_{s}X_s ^N ds+{\displaystyle\sum_{i=1}^{d} K_{s}^{i} X_{s}^N}dW_{s}^{i}, & s \in \left[ t,N\right]  \\
X_t ^N=r^N _t. &
\end{array}
\right.
\end{equation*}
By \textit{duality} between $r^N$ and the process $X^N$, and by estimate (\ref{stimadatkoinf}) we get
\begin{align*}
\Et |r^N _t|^2 & = \Et \int_t ^N \langle \overline{P}_s f_s ,X_s ^N \rangle ds \nonumber  \\
& \leq C \int_t ^N \Vert \overline{P}_s \Vert _{L^{\infty}(\Omega)} \Vert f_s \Vert _{L^{\infty}(\Omega)} e^{-\frac{a}{2}(s-t)}(\Et |r^N _t|^2)^{\frac{1}{2}}  ds \nonumber\\
& \leq \frac{C}{\mu}\int_t ^N e^{-\frac{a}{2}(s-t)} \Vert \overline{P}_s \Vert^2 _{L^{\infty}(\Omega)} \Vert f_s \Vert^2 _{L^{\infty}(\Omega)}ds+ \mu \frac{2}{a} \Et |r^N _t|^2,  \nonumber \\
\end{align*}
where we can take $\mu >0$ such that $\mu \dfrac{2}{a}=\dfrac{1}{2}$. So we get
\begin{equation}\label{stima_r_N}
|r^N _t|^2 \leq C,
\end{equation}
where now $C$ is a constant depending on $a$, $\overline{P}$, $f$,
but $C$ does not depend on $N$. \\ So also $\sup_{t \geq 0}\E |r^N
_t|^2\leq C $. By computing $d|r^N _t|^2$, see e.g. relation
(\ref{Ito^2}) and by the previous estimate we get for every fixed
$T>0$,
\begin{equation*}
 \E \int_0 ^T \sum_{i=1}^{d} |g^{i,N} _s|^2 ds \leq C,
\end{equation*}
where $C>0$ does not depend on $N$. Then we can conclude that for every fixed $T>0$ there exists $\bar r$ and $\bar g$ such that $r^N \rightharpoonup \bar r$ in $L^2_\P(\Omega \times [0,T], \R ^n)$ and $g^N \rightharpoonup \bar g$ in $L^2_\P(\Omega \times [0,T], \R ^{n\times d})$. Moreover, for every $0\leq t\leq T<N$, by the weak convergence,
\begin{equation*}
 \int_t^T H^*_sr^N_s ds \rightarrow \int_t^T H^*_s\bar r_s ds,\qquad 
\int_t^T \sum_{i=1}^d(K^i_s)^*g^{i,N}_s ds \rightarrow \int_t^T \sum_{i=1}^d(K^i_s)^*\bar g^{i}_s ds.
\end{equation*}
By the It\^{o} isometry, the stochastic integral $\int_t^T \sum_{i=1}^d(g^{i,N}_s)^* dW^i_s$ 
is bounded in $L^2_{\F_T}(\Omega)$, so it converges weakly in $L^2_{\F_T}(\Omega)$; 
since every $\Gamma \in L^2_{\F_T}(\Omega)$ can be written as $\Gamma=\int_0^T \gamma^i_s dW^i_s$, weak converges and 
It\^{o} isometry imply that 
\begin{align*}
 &\E\int_0^T \gamma^i_s dW^i_s \int_t^T \sum_{i=1}^d(g^{i,N}_s)^* dW^i_s =\E\int_t^T \sum_{i=1}^d\gamma^i_s ( g^{i,N}_s)^*ds \\
&\rightarrow \E\int_t^T \sum_{i=1}^d\gamma^i_s (\bar g^{i}_s)^*ds=\E\int_0^T \gamma^i_s dW^i_s \int_t^T \sum_{i=1}^d(\bar g^{i}_s)^* dW^i_s,
\end{align*}
and this allows to say that 
\begin{equation*}
 \int_t^T \sum_{i=1}^d(g^{i,N}_s)^* dW^i_s\rightarrow \int_t^T \sum_{i=1}^d(\bar g^{i}_s)^* dW^i_s.
\end{equation*}
So the pair $(\bar r,\bar g)$ is a solution to the elliptic dual equation (\ref{duale}), indeed
\begin{equation*}
\bar r_{t}  =\bar r_{T} + \int_{t}^{T} H_{s} ^{*}\bar r_{s} ds+\int_{t}^{T} \overline{P}_{s}f_{s} ds+ \int_{t}^{T}\sum_{i=1}^{d}  (K_{s}^{i})^{*}\bar g_{s}^{i} ds - \int_{t}^{T} \sum_{i=1}^{d} (\bar g_{s}^{i})^{*} dW_{s}^{i}.
\end{equation*}
Since $T>0$ is arbitrarily, $(\bar r,\bar g)$ is defined on the whole $[0,+ \infty)$. It remains to prove that $\bar{r} \in L^2_\P (\Omega \times [0,+\infty), \R^n)$.
We set
\begin{equation*}
\eta _t ^N =\left\lbrace
\begin{array}[c]{ll}
\bar r_t & 0 \leq t \leq N, \\
0 & t >N.
\end{array}
\right.
\end{equation*}
So $\eta _t ^N \in L^2 _\P(\Omega \times [0,+\infty), \R^n)$. We write a \textit{duality relation}, see (\ref{reldualita}) between $\bar r$ and $X^{\eta ^N}$ solution of the following stochastic differential equation
\begin{equation*}
\left\lbrace
\begin{array}[c]{ll}
 dX_{s}^{\eta ^N }=H_{s}X_{s}^{\eta ^N}ds+\displaystyle\sum_{i=1}^{d}K_{s}^{i}X_{s}^{\eta ^N}dW_{s}^{i}+\eta^N _{s}ds,\\
X_{t}^{\eta^N}=0. & \\
\end{array}
\right.
\end{equation*}
By \textit{duality} we get
\begin{equation*}
\E \int _0 ^N \vert \bar r_s \vert ^2 ds=\E \int_0 ^N \langle \overline{P}_sf_s, X_{s}^{\eta ^N } \rangle ds+\E \langle \bar r_N, X_{N}^{\eta ^N }\rangle.
\end{equation*}
Letting $N\rightarrow\infty$ on both sides we get on the left hand side
\begin{equation*}
\underline{\lim}_{N\rightarrow\infty}\E \int _0 ^N \vert \bar r_s \vert ^2 ds=\lim _{N\rightarrow\infty}\E \int _0 ^N \vert \bar r_s \vert ^2 ds=\E \int _0 ^{+\infty} \vert \bar r_s \vert ^2 ds,
\end{equation*}
by monotone convergence.  On the right hand side, by theorem \ref{remarkDatko} and estimate (\ref{stimanonlin}), we get
\begin{align*}
& \underline{\lim}_{N\rightarrow\infty}\left[\E \int_0 ^N \langle P_sf_s, X_{s}^{\eta ^N } \rangle ds + \E \langle \bar r_N, X_{N}^{\eta ^N }\rangle\right] \nonumber \\
& \leq \underline{\lim} _{N\rightarrow\infty}\frac{C^2}{2}\Vert P \Vert^2_{L^{\infty} (\Omega \times [0,+\infty))}\E\int _0 ^N \vert f_s\vert^2 ds +\frac{1}{2C^2}\E\int _0 ^N \vert X_{s}^{\eta ^N }\vert ^2 ds+\frac{a^2}{2}\E \vert \bar r_N\vert^2+\frac{1}{2a^2}\E \vert X_{N}^{\eta ^N }\vert ^2\nonumber \\
& \leq \frac{C^2}{2}\Vert P \Vert^2_{L^{\infty} (\Omega \times [0,+\infty))} \Vert f\Vert^2 _{L^{2} (\Omega \times [0,+\infty))}ds +C_1\int _0 ^{+\infty} \vert \bar r _s\vert ^2 ds+C+a_1\E \int _0 ^{\infty} \vert \bar r _s \vert^2 ds \nonumber \\
\end{align*}
where we can choose $C_1=a_1=\frac{1}{4}$, and $C$ does not depend on $N$.


Putting together these inequalities we get
\begin{equation*}
 \E \int _0 ^{\infty} \vert \bar r _s \vert^2 ds \leq C^2 \Vert P \Vert^2_{L^{\infty} (\Omega \times [0,+\infty))} \Vert f\Vert^2 _{L^{2} (\Omega \times [0,+\infty))}+2C,
\end{equation*}
and this concludes the proof.
\finedim
\begin{remark}
 As a consequence of the previous proof, we get
\begin{equation}\label{limite.rsegnato}
 \E \vert \bar r _T \vert ^2 \rightarrow 0 \text{ as } T \rightarrow \infty.
\end{equation}
\end{remark}
\begin{remark}
Equation \eqref{duale} has non Lipschitz coefficients and  is a
multidimensional BSDE thus we can not use the Girsanov Theorem, as
done in \cite{BriCon1}, to get rid of the terms involving $K$.
Moreover the typical monotonicity assumptions on the coefficients
of this infinite horizon BSDE, see \cite{BriHu}, are replaced by
the finite cost condition and by the requirement that the minimal
solution $(\overline{P},\overline{Q})$ of \eqref{RiccatiAlg.intro} stabilize
the coefficients relatively the identity, see definition
\ref{def.stab.identita}.
\end{remark}

\section{Synthesis of the optimal control in the infinite horizon case}
We consider the following stochastic differential equation for $t\geq0$:
\begin{equation}
\left\{
\begin{array}
[c]{ll}
dX_{s}=(A_{s}X_{s}+B_{s}u_{s})ds+
{\displaystyle\sum_{i=1}^{d}}
\left(  C_{s}^{i}X_{s}+D_{s}^{i}u_{s}\right)  dW_{s}^{i} +f_s ds & s\geq t\\
X_{t}=x, &
\end{array}
\right.  \label{f.stato.stato}%
\end{equation}
Our purpose is to minimize with respect to $u$ the cost functional,
\begin{equation}
J_{\infty}(0,x,u)=\mathbb{E}\int_{0}^{+\infty}[\left\langle S_{s}X_{s}%
^{0,x,u},X_{s}^{0,x,u}\right\rangle +|u_{s}|^{2}]ds
\label{costo}
\end{equation}
We also introduce the following random variables, for $t\in\lbrack0,+\infty)$:%

\[
J_{\infty}(t,x,u)=\mathbb{E}^{\mathcal{F}_{t}}\int_{t}^{+\infty}[\left\langle
S_{s}X_{s}^{t,x,u},X_{s}^{t,x,u}\right\rangle +|u_{s}|^{2}]ds
\]

Throughout this section and the next section we assume that
\begin{hypothesis}\label{ipotesi_erg}
We will make the following assumptions:
\begin{itemize}
\item $S\geq \epsilon I$, for some $\epsilon >0$. \item
$(A,B,C,D)$ is stabilizable relatively to $S$. \item The first
component $\overline{P}$ of the minimal solution of the Riccati equation is bounded in
time.
\end{itemize}
\end{hypothesis}
Notice that these conditions implies that
$(\overline{P},\overline{Q})$ stabilize $(A,B,C,D)$ relatively to the
identity.
\begin{theorem}
\label{teocontrollo}
Let hypotheses \ref{genhyp}, \ref{hypf} and \ref{ipotesi_erg} hold true, then:

\begin{enumerate}
\item  there exists a unique optimal control $\overline{u}\in L^{2}_\P\left(
\Omega\times\left[  0,+\infty\right)  ,\mathbb{R}^{k}\right)  $ such that
\[
J_{\infty}\left(  0,x,\overline{u}\right)  =\inf_{u\in L^{2}_\P\left(
\Omega\times\left[  0,+\infty\right)  ,\mathbb{R}^{k}\right)  }J_{\infty
}\left(  0,x,u\right)
\]

\item  If $\overline{X}$ is the mild solution of the state
equation corresponding to $\overline{u}$ (that is the optimal
state), then $\overline{X}$ is the unique mild solution to the
\emph{ closed loop } equation for:
\begin{equation}
\left\lbrace
\begin{array}
[c]{ll}
d\overline{X}_{t}=\left[  A_{t}\overline{X}_{t}-B_{t} \left( \Lambda(t, \overline{P}_t , \overline{Q}_t)\overline{X}_{t}+( I+
{\displaystyle\sum_{i=1}^{d}}
\left(  D_{t}^{i}\right)  ^{\ast}\overline{P}_{t}D_{t}^{i})^{-1}(B_{t}^{*}r_{t}+{\displaystyle\sum_{i=1}^{d}}
\left( D_{t}^{i}\right) ^{*}g_{t}^{i})\right)\right]  dt + f_t dt& \\
{\displaystyle\sum_{i=1}^{d}}\left[  C_{s}^{i}\overline{X}_{t}-D_{s}^{i}
(\Lambda(t, \overline{P}_t ,\overline{Q}_t)\overline{X}_{t}+\!\!\!\left(  I+%

{\displaystyle\sum_{i=1}^{d}} \left(  D_{t}^{i}\right)
^{\ast}\overline{P}_{t}D_{t}^{i}\right)
^{-1}\!\!\!\!\!(B_{t}^{*}r_{t}+
{\displaystyle\sum_{i=1}^{d}}\left( D_{t}^{i}\right)
^{*}g_{t}^{i}))\right]   dW^i_{t},\,  t  >0 &
\\
\overline{X}_{0}=x &
\end{array}
\right.
\label{loopf}%
\end{equation}

\item  The following feedback law holds $\mathbb{P}$-a.s. for almost every $t\geq 0$.
\begin{equation}
\overline{u}_{t}=\!-\!\left(  I+ \sum_{i=1}^{d} \left(
D_{t}^{i}\right)  ^{\ast}\overline{P}_{t}D_{t}^{i}\right)
^{-1}\!\!\!\left( \overline{P}_{t}B_{t}+
\sum_{i=1}^{d}\left(\overline{Q}_{t}^{i}D_{t}^{i}+\left(  C_{t}^{i}\right)  ^{\ast}\overline{Q}_{t}D_{t}%
^{i}\right)  \right)  ^{\ast}\!\!\!\overline{X}_{t}+B_{t}^{*}r_{t}+\!\!\sum_{i=1}^{d} (D_{t}^{i})^*g_{t}^{i}. \label{feedback.f}%
\end{equation}

\item  The optimal cost is given by
\begin{align*}
J(0,x,\overline{u}) & =\langle\overline{P}_{0}x,x\rangle+2\langle r_0 ,x\rangle+2\E\int_0 ^{\infty} \langle r_s ,f_s\rangle ds \nonumber \\
& -\E\int_0 ^{\infty} \vert (  I+{\displaystyle\sum_{i=1}^{d}}
\left(  D_{t}^{i}\right)  ^{\ast}\overline{P}_{t}D_{t}^{i})  ^{-1}( B_{t}^{*}r_{t}+{\displaystyle\sum_{i=1}^{d}}\left( D_{t}^{i}\right) ^{*}g_{t}^{i})\vert ^2 ds.
\end{align*}

\end{enumerate}
\end{theorem}
\Dim 
Let us consider the sequence $(\bar{X}^N,\bar{u}^N, \bar{P}^N)$ respectively the optimal state, the optimal control and the first component of the solution to the  Riccati equation \eqref{RiccatiN} defined in section 3 associated with the problem: \\
{  \em minimize over all } $u \in L^2_\P ((0,N)\times \Omega; \R^k):$
\[ J_N(0,x,u) = \E \int_0^N [\langle S_s X_s,X_s \rangle + |u_s|^2] \, ds  \]
where $X$ is the solution to equation \eqref{f.stato.stato} in $[0,N]$.

And let us  consider the minimal solution $\bar P$ evaluated at time $N$ and the associated problem: \\
{  \em minimize over all } $u \in L^2_\P ((0,N)\times \Omega; \R^k):$
\[ J(0,x,u) = \E \int_0^N [\langle S_s X_s,X_s \rangle + |u_s|^2] \, ds + \langle \bar P_NX_N, X_N \rangle .  \]
where $X$ is the solution to equation \eqref{f.stato.stato}.
It turns out from Theorem \ref{controlloT} that the optimal state is the solution to equation \eqref{loopf} considered in $[0,N]$ and that the optimal control is given by \eqref{feedback.f}.
Hence arguing exactly as in Lemma \ref{lemmaSHS} we have that for every $T \leq N$,
\begin{itemize}
\item [a)] $\bar{u}^N \to \bar{u} $ as $ N \to +\infty$ in the space $L^2_\P ((0,T)\times \Omega; \R^k)$ ;
\item [b)] $\bar{X}^N \to \bar{X} $ as $ N \to +\infty$ in the space $L^2_\P ((0,T)\times \Omega; \R^n)$  and then using a) and the Gronwall Lemma,
$\bar{X}^N \to \bar{X} $ as $ N \to +\infty$ in the space $L^2_\P (\Omega;C([0,T];\R^n))$; \end{itemize}

By computing $d\langle {\bar{P}}^N_s \bar{X}^N_s, \bar{X}^N_s\rangle +2\langle  r^N _s,\bar{ X}^N_s\rangle$, where $r^N$ was defined as the first component of the solution 
 $(r^N, g^N)$ of equation \eqref{dualeN} we get,
\begin{align}
& \E\int_0^{N} [\left\langle S_{s}\bar{X}^N_{s},\bar{X}^N_{s}\right\rangle +|\bar{u}^N_{s}|^{2}]ds =
\langle P^N_0 x,x\rangle +2\langle r^N_0 ,x\rangle-2\E \int_0^N \langle  r^N _s,f_s\rangle ds \nonumber \\
&-\E \int_0^N \vert\left(  I+\sum_{i=1}^{d}\left(  D_{s}^{i}\right)^*{P^N}_s D_s^i \right) ^{-1}( B^*_s  r^N_s+\sum_{i=1}^{d}\left(  D_{s}^{i}\right)  ^{\ast} g^{i,N}_{s})\vert  ^2 ds.
\end{align}
and so being $ P^N$ and $r^N$ uniformly bounded, see also \eqref{stima_r_N} and $ f \in L^{\infty}_\P(\Omega \times [0,+\infty), \R^{n })$, we get: 
\begin{align*}
& \E\int_0^{T} [\left\langle S_{s}\bar{X}^N_{s},\bar{X}^N_{s}\right\rangle +|\bar{u}^N_{s}|^{2}]ds \leq
\E\int_0^{N} [\left\langle S_{s}\bar{X}^N_{s},\bar{X}^N_{s}\right\rangle +|\bar{u}^N_{s}|^{2}]ds \leq C 
\end{align*}
and 
\begin{align*}
& \E\int_0^{T} \vert\left(  I+\sum_{i=1}^{d}\left(  D_{s}^{i}\right)^*{P^N}_s D_s^i \right) ^{-1}( B^*_s  r^N_s+\sum_{i=1}^{d}\left(  D_{s}^{i}\right)  ^{\ast} g^{i,N}_{s})\vert  ^2 ds\leq   \\ & 
\E\int_0^{N}  \vert\left(  I+\sum_{i=1}^{d}\left(  D_{s}^{i}\right)^*{P^N}_s D_s^i \right) ^{-1}( B^*_s  r^N_s+\sum_{i=1}^{d}\left(  D_{s}^{i}\right)  ^{\ast} g^{i,N}_{s})\vert  ^2 ds \leq C 
\end{align*}
We get, passing to the limit:
\begin{align*}
& \E\int_0^{T} [\left\langle S_{s}\bar{X}_{s},\bar{X}_{s}\right\rangle +|\bar{u}_{s}|^{2}]ds \leq C
\end{align*}
and letting $ T \to +\infty$:
\begin{align}\label{stima_unif_f}
& \E\int_0^{+\infty} [\left\langle S_{s}\bar{X}_{s},\bar{X}_{s}\right\rangle +|\bar{u}_{s}|^{2}]ds \leq C
\end{align}
Therefore using estimate \eqref{stima_unif_f} and the equation \eqref{f.stato.stato} we get:
\begin{equation}\label{stima_uniforme_X}
\sup_{T >0 }\E |\bar {X}_T|^2 \leq C
\end{equation}
Now we consider the fundamental relation in $[0,T]$:
\begin{align}\label{rel.fond.N}
& \E\int_0^{T} [\left\langle S_{s}\bar{X}^N_{s},\bar{X}^N_{s}\right\rangle +|\bar{u}^N_{s}|^{2}]ds =
\langle P^N_0 x,x\rangle - \langle P^N_T X_T,X_T\rangle  +2\langle r^N_0 ,x\rangle -2\langle r^N_T ,X_T\rangle-2\Et \int_t^T\langle  r^N _s,f_s\rangle ds \nonumber \\
&-\E \int_0^T \vert\left(  I+\sum_{i=1}^{d}\left(  D_{s}^{i}\right)^*{P^N}_s D_s^i \right) ^{-1}( B^*_s  r^N_s+\sum_{i=1}^{d}\left(  D_{s}^{i}\right)  ^{\ast} g^{i,N}_{s})\vert  ^2 ds,
\end{align}
and we notice as first that since all the other terms converge as $N \to +\infty$ we do have,  recalling also that $ r^N \rightharpoonup \bar r$ in $L^2_\P (\Omega \times (0,T); \R^n)$ and $ g^N \rightharpoonup \bar g$ in $L^2_\P (\Omega \times (0,T); \R^k)$,
\begin{align*}
& \lim_{N\to +\infty} \E \int_0^T \vert\left(  I+\sum_{i=1}^{d}\left(  D_{s}^{i}\right)^*{P^N}_s D_s^i \right) ^{-1}( B^*_s  r^N_s+\sum_{i=1}^{d}\left(  D_{s}^{i}\right)  ^{\ast} g^{i,N}_{s})\vert  ^2 ds = \\ &
\E \int_0^T \vert\left(  I+\sum_{i=1}^{d}\left(  D_{s}^{i}\right)^*{\bar{P}}_s D_s^i \right) ^{-1}( B^*_s\bar { r}_s+\sum_{i=1}^{d}\left(  D_{s}^{i}\right)  ^{\ast} \bar{g}^{i}_{s})\vert  ^2 ds
\end{align*}
Hence letting $N$ tend to $+\infty$ in \eqref{rel.fond.N} we get:
\begin{align}\label{rel.fond.Lim}
& \E\int_0^{T} [\left\langle S_{s}\bar{X}_{s},\bar{X}_{s}\right\rangle +|\bar{u}_{s}|^{2}]ds =
\langle \bar{P}_0 x,x\rangle - \langle \bar{P}_T X_T,X_T\rangle  +2\langle \bar{r}_0 ,x\rangle -2\langle\bar{ r}_T ,X_T\rangle-2\Et \int_t^T\langle \bar{ r} _s,f_s\rangle ds \nonumber \\
&-\E \int_0^T \vert\left(  I+\sum_{i=1}^{d}\left(  D_{s}^{i}\right)^*{\bar{P}}_s D_s^i \right) ^{-1}( B^*_s \bar{ r}_s+\sum_{i=1}^{d}\left(  D_{s}^{i}\right)  ^{\ast} \bar{g}^{i}_{s})\vert  ^2 ds,
\end{align}
Therefore, since $\langle \bar{P}_T X_T,X_T\rangle  \geq 0$, and also by \eqref{stima_uniforme_X} and \eqref{limite.rsegnato} we obtain that
\begin{align}\label{prima.dis}
& \E\int_0^{+\infty} [\left\langle S_{s}\bar{X}_{s},\bar{X}_{s}\right\rangle +|\bar{u}_{s}|^{2}]ds \leq 
\langle \bar{P}_0 x,x\rangle +2\langle \bar{r}_0 ,x\rangle -2\E \int_0^{+\infty}\langle \bar{ r} _s,f_s\rangle ds \nonumber \\
&-\E \int_0^{+\infty} \vert\left(  I+\sum_{i=1}^{d}\left(  D_{s}^{i}\right)^*{\bar{P}}_s D_s^i \right) ^{-1}( B^*_s \bar{ r}_s+\sum_{i=1}^{d}\left(  D_{s}^{i}\right)  ^{\ast} \bar{g}^{i}_{s})\vert  ^2 ds,
\end{align}
Now we need to prove the opposite inequality. 
We can reduce to an admissible control such that:
\[ +\infty  > J(0,x,u) = \E \int_0^{+\infty } [|\sqrt{S}_sX^{u,x,0}_s|^2 + |u_s|^2] \, ds   \] 
Hence we have that:
\[  +\infty  > J(0,x,u) = \E \int_0^{+\infty } [|\sqrt{S}_sX^{u,x,0}_s|^2+ |u_s|^2] \, ds \geq 
\E \int_0^{+\infty } [\varepsilon  |X^{u,x,0}_s|^2+ |u_s|^2] \, ds  \]
Thus again the same estimate holds for the state variable $X^{u,x,0}:= X$, 
\begin{equation}\label{stima_uniforme_X_qualunque}
\sup_{T >0 }\E |{X}_T|^2 \leq C
\end{equation}
Computing $d\langle P_s^N X_s, X_s\rangle +2\langle \bar r _s,X_s\rangle$ we get that:
\begin{align}
& \E\int_0^N [\left\langle S_{s}X_{s},X_{s}\right\rangle +|u_{s}|^{2}]ds =\E\langle P_0^N x,x\rangle
+2\langle \bar r_0 ,x\rangle-2\E \int_0^N \langle \bar r _s,f_s\rangle ds \nonumber \\
& +\E \int_0^N \vert \left(  I+\sum_{i=1}^{d}
\left(  D_{s}^{i}\right)  ^{\ast}P_{s}^N D_{s}^{i}\right)  ^{1/2}\left(u_s +( I+\sum_{i=1}^{d}
\left(  D_{s}^{i}\right)  ^{\ast}P_{s}^N D_{s}^{i})  ^{-1}\right.* \nonumber \\
&\left. *\left(
P_{s}^N B_{s}+{\displaystyle\sum_{i=1}^{d}}\left(  Q_{s}^{i,N}D_{s}^{i}+\left(  C_{s}^{i}\right)  ^{\ast}P_{s}^ND_{s}^{i}\right)  \right)  ^{\ast}{X}_{s}
  +B_s ^* \bar r _s +\sum_{i=1}^d D^i_s(\bar g^i_s)^*\right)\vert ^2 ds \nonumber \\
&-\E \int_0^N \vert\left(  I+\sum_{i=1}^{d}\left(  D_{s}^{i}\right)^*P_s^N D_s^i \right) ^{-1}( B^*_s \bar r_s+\sum_{i=1}^{d}
\left(  D_{s}^{i}\right)  ^{\ast}\bar g_{s}^{i})\vert  ^2 ds.
\label{quasirelfond.new}
\end{align}
Now we observe that
\begin{align*}
&\E \int_0^{N}\vert \left(  I+\sum_{i=1}^{d}\left(  D_{s}^{i}\right)^*P_s^N D_s^i \right) ^{-1}( B^*_s \bar r_s+\sum_{i=1}^{d}\left(  D_{s}^{i}\right)  ^{\ast}\bar g_{s}^{i})\vert  ^2 ds \nonumber \\
&=\E \int_0^{+\infty} \vert\left(  I+\sum_{i=1}^{d}\left(  D_{s}^{i}\right)^* P_s^N D_s^i \right) ^{-1}( B^*_s \bar r_s+\sum_{i=1}^{d}\left(  D_{s}^{i}\right)  ^{\ast}\bar g_{s}^{i})\vert  ^2 ds \nonumber \\
 & -\E \int_N^{+\infty}\vert\left(  I+\sum_{i=1}^{d}\left(  D_{s}^{i}\right)^*P_s^N D_s^i \right) ^{-1}(B^*_s \bar r_s+\sum_{i=1}^{d}\left(  D_{s}^{i}\right)  ^{\ast}\bar g_{s}^{i})\vert  ^2 ds,\nonumber \\
\end{align*} 
and by the Dominated Convergence Theorem we have,
\begin{align*}
&\E\int_0^{+\infty} \vert\left(  I+\sum_{i=1}^{d}\left(  D_{s}^{i}\right)^* P_s^N D_s^i \right) ^{-1}( B^*_s \bar r_s+\sum_{i=1}^{d}\left(  D_{s}^{i}\right)  ^{\ast}\bar g_{s}^{i})\vert  ^2 ds \\
&\rightarrow\E\int_0^{+\infty} \vert\left(  I+\sum_{i=1}^{d}\left(  D_{s}^{i}\right)^* \bar{P}_s  D_s^i \right) ^{-1}( B^*_s \bar r_s+\sum_{i=1}^{d}\left(  D_{s}^{i}\right)  ^{\ast}\bar g_{s}^{i})\vert  ^2 ds.\\
\end{align*}
Moreover
\[  \int_N^{+\infty}\vert\left(  I+\sum_{i=1}^{d}\left(  D_{s}^{i}\right)^*P_s^N D_s^i \right) ^{-1}(B^*_s \bar r_s+\sum_{i=1}^{d}\left(  D_{s}^{i}\right)  ^{\ast}\bar g_{s}^{i})\vert  ^2 ds \leq   \int_N^{+\infty}\vert(B^*_s \bar r_s+\sum_{i=1}^{d}\left(  D_{s}^{i}\right)  ^{\ast}\bar g_{s}^{i})\vert  ^2 ds, \]
and, from \eqref{rel.fond.Lim}  we have:
\[  \E\int_0^{+\infty}\vert(B^*_s \bar r_s+\sum_{i=1}^{d}\left(  D_{s}^{i}\right)  ^{\ast}\bar g_{s}^{i})\vert  ^2 ds < +\infty. \]
So we can conclude:
\begin{align*}
 &  \lim_{N\to +\infty} \E \int_N^{+\infty}\vert\left(  I+\sum_{i=1}^{d}\left(  D_{s}^{i}\right)^*P_s^N D_s^i \right) ^{-1}(B^*_s \bar r_s+\sum_{i=1}^{d}\left(  D_{s}^{i}\right)  ^{\ast}\bar g_{s}^{i})\vert  ^2 ds =0\nonumber \\
\end{align*} 
So, letting $N \to +\infty$ in \eqref{quasirelfond.new} we get for every admissible control $u$
\begin{align}\label{seconda.dis}
& \E\int_0^{+\infty} [\left\langle S_{s}X_{s},X_{s}\right\rangle +|u_{s}|^{2}]ds \geq\E\langle \bar P_0 x,x\rangle+2\E\langle \bar r_0 ,x\rangle-2\E \int_t^{+\infty} \langle \bar r _s,f_s\rangle ds \nonumber \\
&-\E \int_0^{+\infty} \vert\left(  I+\sum_{i=1}^{d}\left(  D_{s}^{i}\right)^*\bar P_s D_s^i \right) ^{-1}
( B^*_s \bar r_s+\sum_{i=1}^{d}\left(  D_{s}^{i}\right)  ^{\ast}\bar g_{s}^{i})\vert  ^2 ds.
\end{align}
Then from \eqref{prima.dis} and \eqref{seconda.dis} the theorem easily follows.

\finedim

\section{Ergodic control}
In this section we consider a cost functional depending only on the
asymptotic behaviour of the state (ergodic control). To do it we
first consider a discounted cost functional that fit the assumptions
of section 5 and then we compute a suitable limit of the
discounted cost. Namely, we consider the discounted cost
functional

\begin{equation}
J^{\alpha}(0,x,u)=\mathbb{E}\int_{0}^{+\infty}e^{-2\alpha s}[\left\langle S_{s}X_{s}%
^{0,x,u},X_{s}^{0,x,u}\right\rangle +|u_{s}|^{2}]ds,
\label{costo_scontato}
\end{equation}
where $X$ is solution to equation (\ref{f.stato.stato}), with $A$,
$B$, $C$ and $D$ satisfying hypothesis \ref{genhyp} and $f\in
L^{\infty}_\P(\Omega\times[0,+\infty),\R^n)$. Moreover we recall that we assume hypothesis \ref{ipotesi_erg}. 
When the coefficients are
deterministic the problem has been extensively studied, see e.g.
\cite{Ben1} and \cite{Tess2}.

 Our purpose is to
minimize the discounted cost functional with respect to every
admissible control $u$. We define the set of admissible controls
as
\begin{equation*}
 \mathcal U^{\alpha}=\left\lbrace u\in L^2_\P(\Omega\times[0,+\infty),\R^k):
\mathbb{E}\int_{0}^{+\infty}e^{-2\alpha s}[\left\langle S_{s}X_{s}%
^{0,x,u},X_{s}^{0,x,u}\right\rangle +|u_{s}|^{2}]ds <
+\infty\right\rbrace .
\end{equation*}
Fixed $\alpha >0$, we define $X^{\alpha}_s=e^{-\alpha s}X_s$ and $u^{\alpha}_s=e^{-\alpha s}u_s$: we note that if $u\in \mathcal U _{\alpha}$, then $u^{\alpha}\in \mathcal U$. Moreover we set $A^{\alpha}_s=A_s-\alpha I$ and $f^{\alpha}_s=e^{-\alpha s}f_s$, and $f^{\alpha}\in L^2_\P(\Omega\times[0,+\infty))\cap L^{\infty}_\P(\Omega\times[0,+\infty))$. $X^{\alpha}_s$ is solution to equation
\begin{equation}
\left\{
\begin{array}
[c]{ll}
dX^{\alpha}_s=(A^{\alpha}_sX^{\alpha}_s+B_su^{\alpha}_s)ds+
{\displaystyle\sum_{i=1}^{d}}
\left(  C_s^{i}X^{\alpha}_s+D_s^{i}u^{\alpha}_s\right)  dW_{s}^{i} +f^{\alpha}_s ds & s\geq 0\\
X^{\alpha}_{0}=x, &
\end{array}
\right.  \label{stato.alfa}%
\end{equation}
By the definition of $X^{\alpha}$, we note that if $(A,B,C,D)$ is stabilizable with respect to the identity,
then $(A^{\alpha},B,C,D)$ also is.
We also denote by $(P^{\alpha},Q^{\alpha})$ the solution of the infinite horizon Riccati equation (\ref{RiccatiAlg}),
with $A^{\alpha}$ in the place of $A$. Since, for $0<\alpha<1$, $A^{\alpha}$ is uniformly bounded in $\alpha$,
 also $P^{\alpha}$ is uniformly bounded in $\alpha$.
Now we apply theorem \ref{teocontrollo} to the control problem for the discounted cost $J^{\alpha}$.
Let us denote by $(r^{\alpha},g^{\alpha})$ the solution of the BSDE obtained by equation (\ref{duale}),
where $f$ is replaced with $f^{\alpha}$, and $H$ and $K$ are replaced respectively by $H^{\alpha}$
 and $K^{\alpha}$. $H^{\alpha}$ and $K^{\alpha}$ are defined as in (\ref{notazionifHK}), with $A^{\alpha}$ and
 $P^{\alpha}$ respectively in the place of $A$ and $P$.
\begin{theorem}
Let hypotheses \ref{genhyp} and \ref{hypf} hold true; assume that
$f\in L^{\infty}_\P(\Omega \times [0,+\infty),\R^n)$ then:

\begin{enumerate}
\item  there exists a unique optimal control $\overline{u}^{\alpha}\in \mathcal U _{\alpha}$ such that
\[
J^{\alpha}\left(  0,x,\overline{u}^{\alpha}\right)  =\inf_{u\in \mathcal U _{\alpha}}J^{\alpha}\left(  0,x,u\right)
\]

\item  The following feedback law holds $\mathbb{P}$-a.s. for almost every $t\geq 0$:
\begin{equation}
\overline{u}^{\alpha}_{t}=-\left(  I+ \sum_{i=1}^{d} \left(
D_{t}^{i}\right)  ^{\ast}P^{\alpha}_{t}D_{t}^{i}\right)
^{-1}\!\!\!\!\!\left( P^{\alpha}_{t}B_{t}+
\sum_{i=1}^{d}\left(Q^{\alpha,i}_{t}D_{t}^{i}+\left(  C_{t}^{i}\right)  ^{\ast}P^{\alpha}_{t}D_{t}%
^{i}\right)  \right)
^{\ast}\!\!\!\overline{X}^{\alpha}_{t}+B_{t}^{*}r^{\alpha}_{t}+\sum_{i=1}^{d}D_{t}^{i}(g_{t}^{\alpha,i})^*,
\end{equation}
where $\overline{X}^{\alpha}$ is the optimal state.
\item  The optimal cost $J^{\alpha}\left(  0,x,\overline{u}^{\alpha}\right):=\overline{J}^{\alpha}(x)$ is given by
\begin{align}
\overline{J}^{\alpha}(x) & =\langle P^{\alpha}_{0}x,x\rangle+2\langle r^{\alpha}_0 ,x\rangle+2\E\int_0 ^{\infty} \langle r^{\alpha}_s ,f^{\alpha}_s\rangle ds \nonumber \\
& -\E\int_0 ^{\infty} \vert (  I+
{\displaystyle\sum_{i=1}^{d}}
\left(  D_{t}^{i}\right)  ^{\ast}P^{\alpha}_{t}D_{t}^{i})  ^{-1}( B_{t}^{*}r^{\alpha}_{t}+{\displaystyle\sum_{i=1}^{d}}\left( D_{t}^{i}\right) ^{*}g_{t}^{\alpha,i})\vert ^2 ds. \label{rel.fond.scontata}
\end{align}
\end{enumerate}
\end{theorem}
The optimal cost $\overline{J}^{\alpha}(x)\rightarrow +\infty$ as $\alpha \rightarrow 0$. We want to compute
 $\lim_{\alpha\rightarrow 0}\alpha \overline{J}^{\alpha}(x)$.
In order to do this, we need some convergence results, the first concerning the Riccati equation. To prove this convergence, we note that, by applying the Datko theorem, we are able to prove estimates independent on $\alpha$.
\begin{remark} \label{remarkDatko_alfa}
By the Dakto Theorem, see theorem \ref{remarkDatko}, we can prove an exponential bound for the process $X^{\alpha,t,x}$ which solves the following equation
\begin{equation*}
 \left\lbrace
\begin{array}[c]{ll}
dX^{\alpha,t,x}_s=(H^{\alpha}_s)X^{\alpha,t,x}_sds+{\displaystyle\sum_{i=1}^{d}\left( K^{\alpha,i}_s \right)^{*}X^{\alpha,t,x}_s}dW_{s}^{i}, & s \geq t  \\
X^{\alpha,t,x}_t =x &
\end{array}
\right.
\end{equation*}
We can conclude that there exist $C,a>0$, independent on $\alpha$ such that for every $s\geq t$:
\begin{equation}\label{stimadatkoalfa}
\Et |X^{\alpha}_s| ^2 \leq C e^{-a(s-t)} |x|^2\qquad \mathbb{P}-\text{a.s.}
\end{equation}
\end{remark}
\begin{lemma}\label{lemmaP_alfa}
Assume that hypothesis \ref{genhyp} holds true,
 that $f\in L^{\infty}_\P(\Omega \times [0,+\infty))$. Then $P^{\alpha}_t\rightarrow \overline{P}_t$ as $\alpha\rightarrow 0$
 for all $t \geq 0$, where $\overline{P}$ is the minimal solution of the BRSE.
\end{lemma}
\Dim We can consider the case $t=0$ without loss of generality.

Since $\langle P_0x,x\rangle$, respectively $\langle
P^{\alpha}_0x,x\rangle$, is the optimal cost of the linear
quadratic control problem with state equation given by
(\ref{f.stato.stato}), respectively by (\ref{stato.alfa}), in the
particular case of $f=0$, and cost functional given by
(\ref{costo}), respectively by (\ref{costo_scontato}), we
immediately get that
\begin{equation*}
 P^{\alpha}_0\leq \overline{P}\qquad \text{for all }\alpha>0.
\end{equation*}
 Moreover we get that
\begin{equation*}
 \langle P^{\alpha}_0x,x\rangle=\E \int _0 ^{+\infty} [\langle S_{s}\widehat{X}^{\alpha}(s),\widehat{X}^{\alpha}(s)\rangle
 +|\widehat{u}^{\alpha}(s)|^{2}],
\end{equation*}
where
\begin{equation*}
 \widehat{u}^{\alpha}=-\left(  I+
\sum_{i=1}^{d}
\left(  D_{t}^{i}\right)  ^{\ast}P^{\alpha}_tD_{t}^{i}\right)  ^{-1}\left(
P^{\alpha}_tB_{t}+
\sum_{i=1}^{d}\left(  Q^{\alpha,i}_tD_{t}^{i}+\left(  C_{t}^{i}\right)  ^{\ast}P^{\alpha}_tD_{t}%
^{i}\right)  \right)  ^{\ast}\widehat{X}^{\alpha}_t,
\end{equation*}
and $\widehat{X}^{\alpha}$ is the state corresponding to the
control $\widehat{u}^{\alpha}$. So the pair
$(\widehat{X}^{\alpha},\widehat{u}^{\alpha})$ is bounded in
$L^{2}_\P(\Omega \times [0,+\infty))\times L^{2}_\P(\Omega \times
[0,+\infty))$, so there exists a sequence $\alpha_j \rightarrow 0$
as $j\rightarrow +\infty$ and a pair $(\widehat{X},\widehat{u})$
such that
$(\widehat{X}^{\alpha_j},\widehat{u}^{\alpha_j})\rightharpoonup
(\widehat{X},\widehat{u})$ in $L^{2}_\P(\Omega \times
[0,+\infty))\times L^{2}_\P(\Omega \times [0,+\infty))$. As a
consequence of this convergence, the process $\widehat{X}$ is
solution to equation (\ref{f.stato.stato}), with control
$\widehat{u}$. So we get
\begin{align*}
\langle \overline{P}_0x,x\rangle & \leq \E \int _0 ^{+\infty} [\langle S_{s}\widehat{X}(s),\widehat{X}(s)\rangle +|\widehat{u}(s)|^{2}] \nonumber \\
& \leq \underline{\lim}_{j\rightarrow + \infty}\E \int _0 ^{+\infty} [\langle S_{s}\widehat{X}^{\alpha_j}(s),
\widehat{X}^{\alpha_j}(s)\rangle +|\widehat{u}^{\alpha_j}(s)|^{2}]  =\underline{\lim}_{j\rightarrow + \infty} \E\langle P^{\alpha_j}_0x,x\rangle. \nonumber \\
\end{align*}\finedim

We remark that we can exploit a sort of separation principle,
typical of the linear quadratic case, that allow to estimate
separately the quadratic part from the linear part. Next we want
to prove that, as $\alpha \rightarrow 0$, the optimal pair for the
discounted control problem, that we denote by
$(\widetilde{X}^{\alpha},\widetilde{u}^{\alpha})$ as in the previous
proof, converges to the optimal pair
$(\overline{X},\overline{u})$, defined in theorem
\ref{teocontrollo}.

\begin{lemma}\label{lemmaSHSalfa}
Assume that hypothesis \ref{genhyp} holds true, that $f\in
L^{\infty}_\P(\Omega \times [0,+\infty),\R^k)$. Then, for every $T>0$,
$\widetilde{X}^{\alpha}\rightarrow \overline{X}$ in $L^{2}_\P(\Omega
\times [0,T],\R^n)$ and
$\widetilde{u}^{\alpha}\rightarrow \overline{u}$ in $L^{2}_\P(\Omega
\times [0,T],\R^k)$ as $\alpha\rightarrow 0$.
\end{lemma}
\Dim
We consider the stochastic Hamiltonian system (\ref{SHS}) and the stochastic Hamiltonian system for the discounted problem
\begin{equation}
\left\{
\begin{array}
[c]{l} d\widetilde{X}^{\alpha}_s=[A^{\alpha}_s\widetilde{X}^{\alpha}_s\!-\!B_s(
B_s^{\ast}y_{\alpha}(s)\!\!+ \!\!{\displaystyle\sum_{i=1}^{d}}
\left( D_{s}^{i}\right)  ^{\ast}z^{\alpha,i}_s)]  ds+
{\displaystyle\sum_{i=1}^{d}} [  C_{s}^{i}\widetilde{X}^{\alpha}_s+D_{s}^{i}(
B_{s}^{\ast}y^{\alpha}_s+ {\!\!\!\displaystyle\sum_{k=1}^{d}}
 (D_{s}^{k})^{\ast}z^{\alpha,k}_s ) ] dW_{s}^{i}+f^\alpha_s ds,\\
dy^{\alpha}_s=-[(A^{\alpha}_s)^{\ast}y^{\alpha}_{s}+\!\!\!{\displaystyle\sum_{i=1}^{d}}
\left(  C_{s}^{i}\right)  ^{\ast}z^{\alpha,i}_s+S_{s}X^{\alpha}_s]
ds+ {\displaystyle\sum_{i=1}^{d}}
z^{\alpha,i}_sdW_{s}^{i},\text{ \ \ \ \ \ \ \ \ \ \ \ }t\leq s\leq T,\\
X^{\alpha}_t=x,\\
y^{\alpha}_T=P^{\alpha}_T\tilde{X}^{\alpha}_T,
\end{array}
\right.  \label{SHS_alfa}
\end{equation}
Proceeding as in lemma \ref{lemmaSHS}, we get
\begin{align*}
 &\mathbb{E}^{\mathcal{F}_{t}}\langle y^{\alpha}_{T}-y_{T},\widetilde{X}^{\alpha}_{T}-\overline{X}_{T}\rangle =
 \mathbb{E}^{\mathcal{F}_{t}}\int_{t}^{T}\alpha \langle y^{\alpha}_{s},\overline{X}_{s}\rangle -\alpha \langle y_{s},\widetilde{X}^{\alpha}_{s}\rangle ds+\mathbb{E}^{\mathcal{F}_{t}}\int_{t}^{T} \langle f^{\alpha}_{s}-f_s,y^\alpha_{s}- y_{s}\rangle ds\\
&-\Et \int_t^T \abs{\sqrt{S_s}(\widetilde{X}^{\alpha}_{s}-\overline{X}_{s})}^2ds-\Et \int_{t}^{T} \abs{B_{s}^{\ast}(y^{\alpha}_{s}-y_{s})
+\sum_{i=1}^d(D_s^i)^*(z^{\alpha ,i}_s-z^{i}_s)}^2 ds,\\ 
\end{align*}
that is
\begin{align*}
& \mathbb{E}^{\mathcal{F}_{t}}\langle
P^{\alpha}_{T}(\widetilde{X}^{\alpha}_{T}-\overline{X}_{T}),\widetilde{X}^{\alpha}_{T}-\overline{X}_{T}\rangle +
\Et \!\!\!\langle
(P^{\alpha}_{T}-P_T)\overline{X}_{T},\widetilde{X}^{\alpha}_{T}-\overline{X}_{T}\rangle\!\! =\!\!
\Et\int_{t}^{T}\!\!\!\!\alpha \langle y^{\alpha}_{s}-y_s,\widetilde{X}^{\alpha}_{s}\rangle ds\\
&-\Et\int_{t}^{T}\alpha  \langle y^{\alpha}_{s},\widetilde{X}^{\alpha}_{s}-\overline{X}_s\rangle ds-\Et \int_{t}^{T} \abs{B_{s}^{\ast}(y^{\alpha}_{s}-y_{s})+\sum_{i=1}^d(D_s^i)^*(z^{\alpha,i}_s-z^{i}_s)}^2 ds.\\
 & - \Et \int_t^T
\abs{\sqrt{S_s}(\widetilde{X}^{\alpha}_{s}-\overline{X}_{s})}^2ds
+\mathbb{E}^{\mathcal{F}_{t}}\int_{t}^{T} \langle f^{\alpha}_{s}-f_s,y^\alpha_{s}- y_{s}\rangle ds
\end{align*} It follows
that
\begin{equation*}
 \Et \int_t^T \abs{\sqrt{S_s}(X^{\alpha}_{s}-X_{s})}^2ds+\Et \int_{t}^{T} \abs{B_{s}^{\ast}(y^{\alpha}_{s}-y_{s})+\sum_{i=1}^d(D_s^i)^*(z^{\alpha,i}_s-z^{i}_s)}^2 ds\rightarrow 0
\end{equation*}
as $\alpha \rightarrow 0$.
\finedim

Finally we need to investigate the convergence of $r^{\alpha}$ to $r$, where $(r,g)$ is the solution of equation
(\ref{duale}).
\begin{lemma}
 For all fixed $T>0$, $r^{\alpha}\mid_{[0,T]}\rightarrow r\mid_{[0,T]}$ in $L^2 (\Omega \times [0,T])$.
\end{lemma}
\Dim
First we note that $f^{\alpha}$ is uniformly bounded in $\alpha$ and
\begin{equation}
\E^{\F_{\tau}} \displaystyle\int_{\tau}^{T}\vert H^{\alpha}_{t} \vert^{2}dt+\E^{\F_{\tau}} \displaystyle\int_{\tau}^{T}\vert K^{\alpha}_{t} \vert^{2}dt \leq C ,
\end{equation}
where $C$ is a constant depending on $T$, $x$, $A$, $B$, $C$ and $D$, but not on $\alpha$. So, see proposition \ref{reg-r-infty}, equation (\ref{duale}), where $f$ is replaced by $f^{\alpha}$, and $H$ and $K$ are replaced respectively by $H^{\alpha}$ and $K^{\alpha}$ admits a solution $( r^{\alpha},  g^{\alpha}) \in  L^2_\P (\Omega \times [0,+\infty), \R^n)\times L^2_\P(\Omega \times [0,T], \R^{n \times d}) $, for every $T>0$.
Now let us consider $ \eta \in L^2_\P (\Omega \times [0,T],\R^n)$. $\eta$ can be defined on the whole halfiline  $[0,+\infty)$: we set $\eta_t =0$ for $t>T$. Let $X^{t,x,\eta}$ be the solution of equation (\ref{eqdualita}) and let $X^{\alpha,t,x,\eta}$ be the solution of an equation obtained by equation (\ref{eqdualita}) by replacing $H$ with $H^{\alpha}$ and $K$ with $K^{\alpha}$. By relation (\ref{reldualita}), we get
\begin{equation}
 \E \int_0^T \langle r^{\alpha}_s,\eta_s \rangle ds= \E \int_{0}^{T} \langle P^{\alpha}_sf^{\alpha}_{s}, X_{s}^{\alpha,0,0,\eta}\rangle ds+\E \int_{T}^{+\infty} \langle P^{\alpha}_sf^{\alpha}_{s}, X_{s}^{\alpha,T,X_T^{\alpha,0,0,\eta},0}\rangle ds.\label{alfa_duale}
\end{equation}
and also
\begin{equation}
\E \int_0^T \langle \bar r_s,\eta_s \rangle ds= 
\E \int_{0}^{T} \left\langle \overline{P}_sf_{s}, X_{s}^{0,0,\eta}\right\rangle ds+\E \int_{T }^{+\infty} \left\langle \overline{P}_sf_{s}, X_{s}^{T,X_T^{0,0,\eta},0}\right\rangle ds. \label{limite_duale}
\end{equation}
By theorem \ref{remarkDatko} and by lemmas \ref{lemmaP_alfa} and \ref{lemmaSHSalfa} the right hand side in (\ref{alfa_duale}) converges to the right hand side of (\ref{limite_duale}). So we get that $r^{\alpha}\mid_{[0,T]}\rightharpoonup r\mid_{[0,T]}$ in $L^2_\P (\Omega \times [0,T],\R^n)$. In order to get that $r^{\alpha}\mid_{[0,T]}\rightarrow r\mid_{[0,T]}$ in $L^2_\P (\Omega \times [0,T],\R^n)$, it suffices to prove that $\Vert r^{\alpha}\mid_{[0,T]}\Vert_{L^2_\P (\Omega \times [0,T],\R^n)}\rightarrow \Vert r\mid_{[0,T]}\Vert_{L^2_\P (\Omega \times [0,T],\R^n)}$. We take in (\ref{alfa_duale}) $\eta_t=r^{\alpha}_t$ for $0\leq t\leq T$. 
We get
\begin{equation*}
\E \int_0^T \abs{r^{\alpha}_s}^2 ds =\E \int_{0}^{T} \langle P^{\alpha}_sf^{\alpha}_{s}, X_{s}^{\alpha,0,0,r^{\alpha}}\rangle ds+\E \int_{T}^{+\infty} \langle P^{\alpha}_sf^{\alpha}_{s}, X_{s}^{\alpha,T,X_T^{\alpha,0,0,r^{\alpha}},0}\rangle ds. \\
\end{equation*}
By remark \ref{remarkDatko} and by lemmas \ref{lemmaP_alfa} and \ref{lemmaSHSalfa} the right hand side converges to
\begin{equation*}
 \E \int_{0}^{T} \langle \overline{P}_sf_{s}, X_{s}^{0,0,\bar r,0}\rangle ds+\E \int_{0}^{+\infty} \langle \overline{P}_sf_{s},
  X_{s}^{T,X_T^{0,0,\bar r ,0},0}\rangle ds,\\
\end{equation*}
and this concludes the proof.
\finedim

\begin{remark}\label{remark_r^alfa}
 Following the proof of proposition \ref{reg-r-infty}, it is easy to check that there exists a constant $C>0$, independent on $\alpha$ such that for every $t>0$, $\abs{r^{\alpha}_t}\leq C$.
\end{remark}
We can now study the convergence of $\alpha \overline{J}^{\alpha}$.
\begin{theorem}\label{teo_erg_bis}
Assume that hypothesis \ref{genhyp} holds true, that $f\in
L^{\infty}_\P(\Omega \times [0,+\infty),\R^n)$. Then
\begin{align*}
 &\underline{\lim}_{\alpha\rightarrow 0}\alpha \overline{J}^{\alpha}(x)
=\underline{\lim}_{\alpha\rightarrow 0}2\alpha \E\int_0^{+\infty}
\langle r^{\alpha}_s, f^{\alpha}_s\rangle ds \\
&-\overline{\lim}_{\alpha\rightarrow 0}\alpha\E\int_0 ^{+\infty} \vert (  I+
{\displaystyle\sum_{i=1}^{d}}
\left(  D_{s}^{i}\right)  ^{\ast}P^{\alpha}_{s}D_{s}^{i})  ^{-1}( B_{s}^{*}r^{\alpha}_{s}+{\displaystyle\sum_{i=1}^{d}}\left( D_{s}^{i}\right) ^{*}g_{s}^{\alpha,i})\vert ^2 ds
\end{align*}
\end{theorem}
\Dim
For every $\alpha >0$, by theorem \ref{teocontrollo},
we get 
\begin{align*}
\overline{J}^\alpha(x) & =\langle\overline{P}^\alpha_{0}x,x\rangle+2\langle r^\alpha_0 ,x\rangle+2\E\int_0 ^{\infty} \langle r^\alpha_s ,f^\alpha_s\rangle ds 
 -\E\int_0 ^{\infty} \vert (  I+{\displaystyle\sum_{i=1}^{d}}
\left(  D_{t}^{i}\right)  ^{\ast}\overline{P}^\alpha_{t}D_{t}^{i})  ^{-1}( B_{t}^{*}r^\alpha_{t}+{\displaystyle\sum_{i=1}^{d}}\left( D_{t}^{i}\right) ^{*}g_{t}^{\alpha,i})\vert ^2 ds.
\end{align*}
So, since $P^\alpha$ and $r^\alpha$ are uniformly bounded in alpha, for every $\alpha >0$
\begin{equation*}
\alpha\E\int_0 ^{\infty} \vert (  I+{\displaystyle\sum_{i=1}^{d}}
\left(  D_{t}^{i}\right)  ^{\ast}\overline{P}^\alpha_{t}D_{t}^{i})  ^{-1}( B_{t}^{*}r^\alpha_{t}+{\displaystyle\sum_{i=1}^{d}}\left( D_{t}^{i}\right) ^{*}g_{t}^{\alpha,i})\vert ^2 ds
\leq C,
\end{equation*}
where $C$ is a constant independent on $\alpha$. We can conclude that,
\begin{align*}
 & \underline{\lim}_{\alpha\rightarrow 0}\alpha \overline{J}^{\alpha}(x)=\underline{\lim}_{\alpha\rightarrow 0}2\alpha \E\int_0^{+\infty}\langle r^{\alpha}_s, f^{\alpha}_s\rangle ds \\
&-\overline{\lim}_{\alpha\rightarrow 0}\alpha\E\int_0 ^{+\infty} \vert (  I+
{\displaystyle\sum_{i=1}^{d}}
\left(  D_{s}^{i}\right)  ^{\ast}P^{\alpha}_{s}D_{s}^{i})  ^{-1}( B_{s}^{*}r^{\alpha}_{s}+{\displaystyle\sum_{i=1}^{d}}\left( D_{s}^{i}\right) ^{*}g_{s}^{\alpha,i})\vert ^2 ds
\end{align*}

\finedim


\begin{thebibliography}{99}

\bibitem{Ben}
A.~Bensoussan.
\newblock {Lectures on Stochastic Control}.
\newblock {in {\em Nonlinear Filtering and Stochastic Control }Proceedings, Cortona
1981}.Lecture Notes in Math., 972, Springer, Berlin 1982

\bibitem{Ben1}
A.~Bensoussan and J.~Frehse.
\newblock {On Bellman Equations of Ergodic Control in $\mathbb{R}^b$}.
\newblock {{\em J. reine angew. Math.} 429(1992),125-160.}


\bibitem{Bi76}
J.-M.~Bismut.
\newblock {Linear quadratic optimal stochastic
control with random coefficients}.
\newblock {\em SIAM J. Contr. Optim.} 14
(1976), 419--444.

\bibitem{Bi78}
J.-M.~Bismut.
\newblock {Contr\^{o}le des syst\`{e}mes lin\'{e}aires quadratiques: applications de l'int\'{e}grale
stochastique.}
\newblock {In {\em S\'{e}minaire de Probabilit\'{e}s, XII}}.
Lecture Notes in Math., 649, Springer, Berlin, 1978.


\bibitem{BriCon}
P.~Briand and F.~Confortola
\newblock {BSDEs with stochastic Lipschitz condition and quadratic PDEs in Hilbert spaces}
\newblock {\em to appear on Stochastic Process. Appl.}

\bibitem{BriCon1}
P.~Briand and F.~Confortola
\newblock {Quadratic BSDEs with random terminal time and elliptic PDEs in infinite dimension}
\newblock {\em  arXiv:0704.1223.}

\bibitem{BriHu}
P.~Briand and Y.~Hu
\newblock {Stability of BSDEs with random terminal time and homogenization of semilinear elliptic
PDEs.}
\newblock {{\em J. Funct. Anal.} 155 (1998), no. 2, 455--494.}

\bibitem{BucPen}
 R. Buckdahn and S. Peng.
\newblock{Stationary backward stochastic differential equations and
associated partial differential equations,}
\newblock{\em Probability Theory and Related Fields,} 115 (1999), pp. 383-399.



\bibitem{DaPIch}
G.~Da Prato and A.~Ichikawa.
\newblock {Quadratic optimal for linear periodic systems}.
\newblock {\em Appl. Math. Optim.} 18
(1988), 39--66.

\bibitem{Datko}
R.~Datko.
\newblock{Extending a Theorem of A. M. Liapunov to Hilbert space,}
\newblock{\em J. Math. Analysis Applic.,} 32 (1970), pp. 610-616.


\bibitem{Gal}
L. I. Gal'chuk.
\newblock{Existence and uniqueness of a solution for stochastic equations with respect to semimartingales}
\newblock{\em Theory of Probability and its Applications} Vol XXII, n. 4, (1978) 751-763.

\bibitem{GTinf}
G.~Guatteri and G. Tessitore.
\newblock{Backward Stochastic Riccati Equations and \\
Infinite Horizon
 L-Q Optimal Control Problems with Stochastic Coefficients}
\newblock{{\em Applied Mathematics and Optimization,}} 57 (2008), no.2, pp.207-235  


\bibitem{Ich}
A.~Ichikawa.
\newblock{Equivalence of $L_p$ stability and exponential stability for a class
 of nonlinear semigroups,}
\newblock{\em Nonlinear Analysis, Theory, Methods \& Applications,} 8 (1984), pp.
805-815.

\bibitem{KoTa01}
M.~Kohlmann and S.~Tang.
\newblock {  New developments in backward stochastic Riccati
 equations and their applications}.
\newblock {In {\em  Mathematical finance (Konstanz, 2000)}},
Trends Math., Birkh�ser, Basel, 2001.



\bibitem{KoTa02}
M.~Kohlmann and S.~Tang.
\newblock { Global adapted solution of
one-dimensional backward stochastic Riccati equations, with
application to the mean-variance hedging}.
\newblock {\em  Stochastic Process. Appl.} 97
(2002),  1255--288.


\bibitem{KoTa03}
M.~Kohlmann and S.~Tang.
\newblock { Multidimensional backward stochastic Riccati equations and
applications}.
\newblock {\em  SIAM J. Contr. Optim.} 41 (2003),  1696--1721.

\bibitem{KoZh}
M.~Kohlmann and X.Y.~Zhou.
\newblock {Relationship between backward stochastic differential equations
and stochastic controls: a linear-quadratic approach. }
\newblock {\em  SIAM J. Contr. Optim.} 38 (2000),  1392--1407.



\bibitem{LeSM}
J.-P.~Lepeltier, J.~San Mart{\'{\i}}n.
   \newblock { Existence for {BSDE} with superlinear-quadratic
   coefficient}.
   \newblock {\em Stochastics Stochastic Rep.},
 {63} (1998),227--240.




\bibitem{Peng}
S.~Peng.
\newblock { Stochastic Hamilton-Jacobi-Bellman  Equations.}
\newblock {\em SIAM J. Contr. Optim.} 30 (1992),  284--304.

\bibitem{Pe99}
S.~Peng. \newblock {Open problems on backward stochastic
differential equations.}
\newblock { In {\em Control of distributed parameter and
stochastic systems (Hangzhou, 1998)}},  Kluwer Acad. Publ.,
Boston, 1999.

\bibitem{PeWu}
S.~Peng, Z.~Wu. \newblock {Fully coupled forward-backward stochastic differential equations and applications to optimal control.}
\newblock { SIAM J. Control Optim. 37 (1999), no. 3, 825--843.}.

\bibitem{Tang}
S.~Tang.
\newblock { General linear quadratic optimal control problems with random coefficients:
linear stochastic Hamilton systems and backward stochastic Riccati
equations.}
\newblock {SIAM J. Control Optim. 42 (2003), no. 1, 53--75}


\bibitem{Tess}
G.Tessitore.
\newblock { Some remarks on the Riccati equation arising in an optimal control
problem with state- and control-dependent noise.}
\newblock {\em SIAM J. Contr. and Optim. 30 (1992),  717--744.}


\bibitem{Tess2}
G.Tessitore.
\newblock { Infinite Horizon, Ergodic and Periodic Control for a Stochastic Infinite Dimensional Affine Equation.}
\newblock {\em Journal of Mathematical Systems, Estimation, and Control. 8 n.4 (1998), 1-28}

\bibitem{Yo-Z}
J. Yong and X.Z. Zhou
\newblock { Stochastic controls. Hamiltonian systems and HJB equations. }
\newblock {\em Applications of Mathematics (New York), 43. Springer-Verlag, New York, 1999}


\end{thebibliography}
\end{document}